\documentclass[10pt]{siamltex}

\usepackage{amssymb,color,amsmath,latexsym,eucal,cite}
\usepackage{a4wide,graphicx,psfrag} 
\usepackage[latin1]{inputenc}
\usepackage{cases}
\usepackage{epstopdf}

\usepackage{algorithm}
\newtheorem{ipotesi}{Assumption}[section]    

\newtheorem{thm}{Theorem}[section]

\newtheorem{lem}[thm]{Lemma}

\newcommand{\R}{\mathbb{R}}

\newcommand{\Ek}{\mathbb{E}_k}
\newcommand{\E}{\mathbb{E}}

\renewcommand\refname{References}

\title{An investigation of  stochastic trust-region based algorithms for finite-sum minimization}
\author{ Stefania Bellavia\footnotemark[1], Benedetta Morini\footnotemark[1], Simone Rebegoldi\footnotemark[2]}

\begin{document}
\maketitle
\footnotetext[1]{Dipartimento  di Ingegneria Industriale, Universit\`a degli Studi di Firenze,
Viale G.B. Morgagni 40,  50134 Firenze,  Italia. Members 
of the Gruppo Nazionale Calcolo Scientifico-Istituto Nazionale di Alta Matematica
(GNCS-INDAM). Emails:
stefania.bellavia@unifi.it, benedetta.morini@unifi.it}

\footnotetext[2]{Dipartimento di Scienze Fisiche, Informatiche e Matematiche, Universit\`a degli Studi di Modena e Reggio Emilia,
Via Campi 213/B, 41125 Modena, Italia. Member 
of the Gruppo Nazionale Calcolo Scientifico-Istituto Nazionale di Alta Matematica
(GNCS-INDAM). Email: simone.rebegoldi@unimore.it}

\footnotetext[3]{ 
The research that led to the present paper was partially supported by INDAM-GNCS through Progetti di Ricerca 2023 and 
by PNRR - Missione 4 Istruzione e Ricerca - Componente C2 Investimento 1.1, Fondo per il Programma Nazionale di Ricerca e Progetti di Rilevante Interesse Nazionale (PRIN) funded by the European Commission under the NextGeneration EU programme, project ``Advanced optimization METhods for automated central veIn Sign detection in multiple sclerosis from magneTic resonAnce imaging (AMETISTA)'',  code: P2022J9SNP,
MUR D.D. financing decree n. 1379 of 1st September 2023 (CUP E53D23017980001), project 
``Numerical Optimization with Adaptive Accuracy and Applications to Machine Learning'',  code: 2022N3ZNAX
 MUR D.D. financing decree n. 973 of 30th June 2023 (CUP B53D23012670006), and project ``Inverse problems in PDE: theoretical and numerical analysis'', code: 2022B32J5C MUR D.D. financing decree n. 973 of 30th June 2023 (CUP B53D23009200006). 
}
\begin{abstract}
This work elaborates on the TRust-region-ish (TRish) algorithm, a stochastic optimization method for finite-sum minimization problems proposed by Curtis et al. in \cite{Curtis2019, Curtis2022}. A theoretical analysis that complements the results in the literature is presented, and the  
issue of tuning the involved  hyper-parameters   is investigated. Our study also focuses on a practical version of the method, which computes the
stochastic gradient by means of the inner product test and the orthogonality 
test proposed by Bollapragada et al. in \cite{Bollapragada2018}. It is shown experimentally that this implementation improves the performance of TRish and reduces its  sensitivity to the choice of the hyper-parameters.
\end{abstract}

\section{Introduction}
We consider the finite-sum minimization problem
\begin{equation}\label{eq:problem}
\min_{x\in\mathbb{R}^n}F(x)\equiv \frac{1}{N}\sum_{i=1}^{N} F_i(x),
\end{equation}
where $N$ is a large integer and each function $F_i$, $1\le i\le N$,  is differentiable, and we address
its solution by the stochastic first-order method named TRust-regionish-algorithm (TRish) and presented in \cite{Curtis2019}.

TRish combines ideas from the Stochastic Gradient (SG) method \cite{RM} and  the trust-region methodology \cite{cgt}.
Depending on the value of the norm of the stochastic gradient, the step takes  either the form 
of the SG step or the form of a normalized step inspired by a first-order trust-region scheme.
This strategy aims at reducing the dependence of the steplength from the specific function $F$, and as a result, TRish can outperform a traditional SG algorithm from the numerical viewpoint, see the experiments in \cite{Curtis2019}. An extension of TRish 
to second-order models was recently proposed in \cite{Curtis2022}.

TRish differs from the stochastic trust-region methods in \cite{BSV,BGMT, Bellavia-et-al-20, Bellavia-et-al-22, bcms,
Cartis-et-al, ChenMeniSche18, CS} in two main respects.
First, the computed step is accepted at each iteration, namely it is not ruled by the actual-to-predicted  
reduction ratio.
Second, the theoretical properties of TRish are obtained by assuming an upper bound on the variance  of the stochastic gradient estimates, which is  common in the analysis of the SG method, whereas specific accuracy requirements on the random model at each iteration are not required.

Our interest in TRish derives from the experimental comparison we carried out among  various stochastic trust-region methods. 
In our experience, after tuning the parameters,  TRish  compares well with  adaptive stochastic trust-region methods, see \cite{Bellavia-et-al-22}. Limiting our attention to the first-order version of TRish, we observe that the computation of the step at each iteration requires three hyper-parameters, as well as a rule for computing the stochastic gradient.
These hyper-parameters are the steplength and two positive scalars employed 
in the choice between  a SG-type step and a normalized step; these parameters can vary along the iterations.
Due to the form of problem \eqref{eq:problem},
the stochastic gradient can be naturally built using a mini-batch approach, which requires the selection of a sample size at each iteration. Note that the sample size is implicitly a further hyper-parameter of the algorithm.
Clearly, the performance of TRish depends on the hyper-parameters choice, and its application  requires tuning them for each problem, analogously to SG methods. This tuning process can be computationally onerous, and represents the main limitation of the algorithm.

In this paper, we pursue the analysis of TRish, by focusing on its theoretical properties and the practical choice of both the hyper-parameters and  the batch size used for  computing
the gradient estimator. We study the theoretical properties of the first-order version of the algorithm, under the assumption that the three hyper-parameters above-mentioned are constant, and provide results that complement the convergence analysis in \cite{Curtis2019, Curtis2022}, in terms of both conditions on such parameters and expected decrease of the optimality measure.
Furthermore, inspired by \cite{Bollapragada2018}, we investigate the use of adaptive sampling techniques
that vary dynamically the sample size for computing the stochastic gradient, and provide a practical implementation of TRish that complies with the theoretical assumptions of our analysis. We show experimentally that an adaptive strategy for building the stochastic gradient reduces the dependence of the performance of TRish from its parameters, although it does not completely overcome the need to tune them. 

The paper is structured as follows. In section \ref{sec2}, we briefly recall TRish and explore its convergence properties under different settings. In section \ref{sec3}, we present an implementation of the algorithm, named TRish\_AS, which is based on the adaptive sampling techniques proposed in \cite{Bollapragada2018}. In section \ref{sec4}, we report the numerical performance of TRish\_AS on a series of classification and regression problems. Some conclusions are offered in section \ref{sec:conclusions}.

{\bf Notation.} The symbol $\Ek[\cdot]$ denotes the conditional 
expectation of a random variable conditioned on the event that the $k-$th iterate of the algorithm is $x_k$.
The symbol $\|\cdot\|$ denotes the Euclidean norm.

\section{The algorithm and convergence results}\label{sec2}
TRish is a simple  yet effective  algorithm for  stochastic and finite-sum minimization problems  \cite{Curtis2019}. At each iteration $k$,  a stochastic estimate $g_k$ of the gradient is formed and used for building 
the step. A normalized step along $g_k$ is built whenever the  norm of $g_k$ lies within a specified interval $\left [\frac{1}{\gamma_{1,k}}, \frac{1}{\gamma_{2,k}}\right]$, otherwise  the step is of SG-type. The stepsize depends on a given positive scalar $\alpha_k$. In  \cite{Curtis2019, Curtis2022}, convergence results have been given for constant sequences $\{\alpha_k\}$, $\{\gamma_{1,k}\}$, $\{\gamma_{2,k}\}$, as well as
for diminishing stepsizes. Here, we focus on the case where the parameters sequences are constant.
The TRish algorithm is given below. 

\begin{algorithm}[h!]
\caption{{\bf TRish}}\label{algo:trish}
\vspace{2mm}
\noindent Choose an initial iterate $x_0\in\mathbb{R}^n$.\\
Choose a steplength $\alpha>0$ and parameters $0<\gamma_2<\gamma_1$.\\
Set $k=0$.
\vskip 2pt\noindent
Repeat until a convergence test is satisfied
\begin{enumerate}
\item  Generate a stochastic gradient $g_k\approx \nabla F(x_k)$.
\item Compute
\begin{equation}\label{eq:trish_step}
p_k=-\begin{cases}
\gamma_1 \alpha g_k, \quad &\text{if }\|g_k\|\in[0,\frac{1}{\gamma_1})\\
\alpha \displaystyle \frac{g_k}{\|g_k\|}, \quad &\text{if }\|g_k\|\in[\frac{1}{\gamma_1},\frac{1}{\gamma_2}] \\
\gamma_2 \alpha g_k, \quad &\text{if }\|g_k\|\in(\frac{1}{\gamma_2},+\infty[
\end{cases}, \qquad x_{k+1}=x_k+p_k.
\end{equation} 
\item Set $k=k+1$.
\end{enumerate}
\end{algorithm}
Algorithm \ref{algo:trish} is based on a {\em trust-regionish} procedure. The step $p_k$ takes different forms depending on the three cases in 
(\ref{eq:trish_step}). In the second case,  $\|g_k\|\in[\frac{1}{\gamma_1},\frac{1}{\gamma_2}] $,
the step taken solves the trust-region problem 
$$
\min_{\begin{array}{l} 
p\in \R^n\\
\|p\|\le \alpha
\end{array}} F(x_k)+g_k^T p.
$$
In contrast to  deterministic and stochastic trust-region methodologies\cite{BSV,BGMT, Bellavia-et-al-20, Bellavia-et-al-22, bcms,
cgt, Cartis-et-al, ChenMeniSche18, CS}, 
the adopted step may differ from the  trust-region step and it  is accepted at each iteration, namely no acceptance rules are tested.
The reason for choosing between a trust-region step and a SG-type step derives from the observation that using the trust-region step
$p_k= -\alpha g_k/\|g_k\|$ at each iteration can  prevent the algorithm  from  making progress in expectation, see \cite[p. 202]{Curtis2019}. 
\vskip 5pt \noindent
In the following, we provide convergence results for TRish that complement 
those given in \cite{Curtis2019, Curtis2022}. We make the following assumptions on $F$.
\vskip 5pt \noindent
\begin{ipotesi}\label{H1}
$F$ is continuously differentiable and bounded from below by $F_*=\inf_{x\in \mathbb{R}^n} F(x)$. The gradient 
$\nabla F:{\mathbb{R}^n}\rightarrow\mathbb{R}^n$ is Lipschitz continuous with constant $L> 0$, i.e., 
\begin{equation}\label{eq:Lipschitz}
\|\nabla F(x)-\nabla F(y)\|\leq L \|y-x\|, \quad \forall \ x,y \in\mathbb{R}^n.
\end{equation}
\end{ipotesi}
\vskip 2pt \noindent
Results are presented for both general nonconvex objective functions $F$ and for  
functions satisfying the Polyak-Lojasiewicz (PL) condition stated below.
\begin{ipotesi}\label{H2}
There exists a constant $\mu>0$ such that 
\begin{equation}\label{eq:PL}
\|\nabla F(x)\|^2\geq 2\mu(F(x)-F_*), \quad \forall \ x\in\mathbb{R}^n.
\end{equation}
\end{ipotesi}
\vskip 5pt
\noindent
Strongly convex functions and strongly convex functions composed with linear maps satisfy the PL condition \cite{Karimi}. We note that $L\geq \mu$, see e.g. \cite[pag. 246]{nocedalsurvey}.

As one would expect, the analysis of TRish heavily relies on the assumptions made on the stochastic gradient $g_k$. In the remainder of this section, we make the following assumption.
\vskip 5pt \noindent
\begin{ipotesi}\label{H3}
The vector $g_k$ is an unbiased estimator of $\nabla F(x_k)$, i.e.,
\begin{equation}\label{eq:unbiased}
\Ek[g_k]=\nabla F(x_k).
\end{equation}
\end{ipotesi}
Under this assumption,  $(- g_k)$ is a descent direction in expectation since 
\begin{equation}\label{discesa}
\Ek[g_k^T\nabla F(x_k)]=\|\nabla F(x_k)\|^2>0,
\end{equation}
and consequently 
\begin{equation}\label{eq:fond4}
\Ek[\|\nabla F(x_k)-g_k\|^2]=-\|\nabla F(x_k)\|^2+\Ek[\|g_k\|^2].
\end{equation}

Before presenting our analysis, we give two simple inequalities that will be repeatedly used 
\begin{eqnarray}
& &0\le \frac 1 2\left(\|g_k\|-\|\nabla F(x_k)-g_k\|\right)^2=\frac 1 2 \|g_k\|^2-\|\nabla F(x_k)-g_k\|\|g_k\|+
\nonumber \\
& & \qquad \frac 1 2 \|\nabla F(x_k)-g_k\|^2, \label{eq:tec1} \\
& & 
0\le \frac{\gamma_2}{2\gamma_1^2} \left(1- \frac{\gamma_1^2}{\gamma_2}  \|\nabla F(x_k)-g_k\|\right)^2=
\frac{\gamma_2}{2\gamma_1^2}-\|\nabla F(x_k)-g_k\| +\nonumber \\
& & \qquad \frac{\gamma_1^2}{2\gamma_2}  \|\nabla F(x_k)-g_k\|^2, \label{eq:tec2}
\end{eqnarray}
and  we introduce a positive scalar that will appear in the results
\begin{equation}\label{beta}
\beta= \frac{\gamma_1^2-\gamma_2^2}{2\gamma_2}+\frac{1}{2} \alpha\gamma_1^2L .
\end{equation}
Note that $\beta$ is positive as $\gamma_2< \gamma_1$.

The first result provides bounds on the expected decrease of the objective function yielded by TRish from one iteration to the other.  

\vskip 5pt
\begin{lem}\label{lem:tec1}
Let Assumptions \ref{H1} and \ref{H3} hold. 
For all $k$, it holds
\begin{eqnarray}
\Ek[F(x_{k+1})]&\le& F(x_k)-\alpha \frac{\gamma_1^2}{2\gamma_2} \|\nabla F(x_k)\|^2+\alpha \beta \Ek[\|g_k\|^2] \label{eq:fond3},\\
\Ek[F(x_{k+1})]&\le& F(x_k)-\frac 1 2 \alpha \left(\gamma_2-\alpha \gamma_1^2 L \right )\|\nabla F(x_k)\|^2+
\nonumber \\
& & \alpha \beta \Ek[\|\nabla F(x_k)-g_k\|^2], \label{eq:fond5}
\end{eqnarray}
with $\beta$ given in (\ref{beta}).
\end{lem}

{\bf Proof.}
By applying the Descent Lemma for functions with a Lipschitz continuous gradient \cite[Proposition A.24]{Bertsekas-2016}, we obtain
\begin{equation}\label{eq:start}
F(x_{k+1})\leq F(x_k)+\nabla F(x_k)^Tp_k+\frac{L}{2}\|p_k\|^2,
\end{equation}
and it trivially follows
\begin{equation}\label{eq:fond1}
F(x_{k+1})\le F(x_k)+g_k^Tp_k+(\nabla F(x_k)-g_k)^T p_k +\frac 1 2 L \|p_k\|^2.
\end{equation}
Then, we show that, for all $k$, it holds
\begin{equation}\label{eq:fond2}
F(x_{k+1})\le F(x_k)-\frac 1 2 \alpha (\gamma_2- \alpha \gamma_1^2 L) \|g_k\|^2  +\alpha\frac{\gamma_1^2}{2\gamma_2} \|\nabla F(x_k)-g_k\|^2. 
\end{equation}
To this aim, we prove that inequality \eqref{eq:fond2} holds irrespectively of the form of $p_k$, namely it holds for Case 1 ($p_k=- \gamma_1\alpha  g_k$), Case 2 ($p_k=-\alpha g_k/\|g_k\|$), and Case 3 ($p_k=-\gamma_2 \alpha g_k$).

In Case 1, we have $p_k=-\gamma_1 \alpha  g_k$
and equations (\ref{eq:fond1}), (\ref{eq:tec1})  give
\begin{eqnarray*}
F(x_{k+1})&\le& F(x_k)-\alpha \gamma_1\|g_k\|^2-\alpha \gamma_1(\nabla F(x_k)-g_k)^Tg_k +\frac 1 2 \alpha^2 \gamma_1^2 L \|g_k\|^2\\
& \le & F(x_k)-\alpha \gamma_1\|g_k\|^2+\alpha \gamma_1\|\nabla F(x_k)-g_k\|\|g_k\| +\frac 1 2 \alpha^2 \gamma_1^2 L \|g_k\|^2\\
&\le & F(x_k)-\alpha \gamma_1\|g_k\|^2+\frac 1 2 \alpha \gamma_1\|g_k\|^2  + \frac 1 2 \alpha\gamma_1 \|\nabla F(x_k)-g_k\|^2 +\frac 1 2 \alpha^2 \gamma_1^2 L \|g_k\|^2\\
&=&  F(x_k)-\frac 1 2 \alpha \gamma_1\|g_k\|^2 + \frac 1 2 \alpha\gamma_1 \|\nabla F(x_k)-g_k\|^2 +\frac 1 2 \alpha^2 \gamma_1^2 L \|g_k\|^2\\
&\le&  F(x_k)-\frac 1 2 \alpha \gamma_2\|g_k\|^2 +  \alpha\frac{\gamma_1^2}{2\gamma_2} \|\nabla F(x_k)-g_k\|^2 +\frac 1 2 \alpha^2 \gamma_1^2 L \|g_k\|^2,
\end{eqnarray*}
where the last inequality follows from $-\gamma_1< -\gamma_2$ and $\gamma_1/\gamma_2>1$.

In Case 2, we have $p_k=-\alpha  g_k/\|g_k\|$
and equations  (\ref{eq:fond1}),  (\ref{eq:tec2})  give
\begin{eqnarray*}
F(x_{k+1})&\le& F(x_k)-\alpha \|g_k\|-\alpha (\nabla F(x_k)-g_k)^T\frac{g_k}{\|g_k\|} +\frac 1 2 \alpha^2  L\\
& \le &F(x_k)-\alpha \|g_k\|+\alpha \|\nabla F(x_k)-g_k\| +\frac 1 2 \alpha^2  L \\
& \le & F(x_k)-\alpha\|g_k\|+\alpha \frac{\gamma_2}{2\gamma_1^2}+ \alpha \frac{\gamma_1^2}{2\gamma_2} \|\nabla F(x_k)-g_k\| ^2
+\frac 1 2 \alpha^2  L.
\end{eqnarray*}
By noting that  $\gamma_1\|g_k\|\ge 1$ and $\gamma_2\|g_k\|\le 1$, it follows
\begin{eqnarray*}
F(x_{k+1})&\le& F(x_k)-\alpha \gamma_2\|g_k\|^2+\alpha \frac{\gamma_2}{2\gamma_1^2}\gamma_1^2 \|g_k\|^2 +
\alpha \frac{\gamma_1^2}{2\gamma_2} \|\nabla F(x_k)-g_k\| ^2+ \frac 1 2 \alpha^2 \gamma_1^2 L\|g_k\|^2\\
& \le & F(x_k)-\frac{1}{2} \alpha\gamma_2 \|g_k\|^2 +  \alpha\frac{\gamma_1^2}{2\gamma_2}  \|\nabla F(x_k)-g_k\|^2 
+\frac 1 2 \alpha^2 \gamma_1^2 L\|g_k\|^2.
\end{eqnarray*}

In Case 3, we have $p_k=-\gamma_2 \alpha g_k$,
and  proceeding as in Case 1,  equations (\ref{eq:fond1}), (\ref{eq:tec1})  give
\begin{eqnarray*}
F(x_{k+1})
&\le & F(x_k)-\frac 1 2 \alpha \gamma_2\|g_k\|^2 + \frac 1 2 \alpha\gamma_2 \|\nabla F(x_k)-g_k\|^2 +
\frac 1 2 \alpha^2 \gamma_2^2 L \|g_k\|^2\\
&\le & F(x_k)-\frac 1 2 \alpha \gamma_2\|g_k\|^2 +  \alpha\frac{\gamma_1^2}{2\gamma_2} \|\nabla F(x_k)-g_k\|^2 +
\frac 1 2 \alpha^2 \gamma_1^2 L \|g_k\|^2,\\
\end{eqnarray*}
where the last inequality follows from $\gamma_2<  \gamma_1$ and $\gamma_1/\gamma_2>1$.

Finally, equations (\ref{eq:fond2}) and (\ref{eq:fond4}) yield
\begin{eqnarray*}
\Ek[F(x_{k+1})]&\le& F(x_k)-  \frac{1}{2}\alpha(\gamma_2-\alpha\gamma_1^2 L) \E[|g_k\|^2]+ 
\alpha \frac{\gamma_1^2}{2\gamma_2}\Ek[\|\nabla F(x_k)-g_k\|^2]\\
&=&  F(x_k)-   \frac{1}{2}\alpha(\gamma_2-\alpha\gamma_1^2 L) \E[|g_k\|^2]+ 
\alpha \frac{\gamma_1^2}{2\gamma_2}(-\|\nabla F(x_k)\|^2+ \Ek[\|g_k\|^2])\\
&=& F(x_k)-\alpha\frac{\gamma_1^2}{2\gamma_2}\|\nabla F(x_k)\|^2+\alpha
\left(\frac{\gamma_1^2-\gamma_2^2}{2\gamma_2}+\frac 1 2\alpha\gamma_1^2L\right)\Ek[\|g_k\|^2],
\end{eqnarray*}
which gives (\ref{eq:fond3}) due to the definition (\ref{beta}) of $\beta$. Furthermore, equations (\ref{eq:fond2}) and (\ref{eq:fond4})  give
\begin{eqnarray*}
\Ek[F(x_{k+1})]&\le& F(x_k)-  \frac{1}{2}\alpha(\gamma_2-\alpha\gamma_1^2 L) \E[\|g_k\|^2]+ 
\alpha \frac{\gamma_1^2}{2\gamma_2}\Ek[\|\nabla F(x_k)-g_k\|^2]\\
&=&  F(x_k)- \frac{1}{2}\alpha(\gamma_2-\alpha\gamma_1^2 L) (\|\nabla F(x_k)\|^2+\Ek[\|\nabla F(x_k)-g_k\|^2]) + \\
& &  \, \alpha \frac{\gamma_1^2}{2\gamma_2}\Ek[\|\nabla F(x_k)-g_k\|^2],
\end{eqnarray*}
and (\ref{eq:fond5}) follows from the definition (\ref{beta}) of $\beta$.  $\Box$
\vskip 5pt
With respect to the existing literature, the arguments used   above  
are simpler than those used in \cite[Lemma 1]{Curtis2019} for 
deriving an analogous result. Furthermore, comparing   (\ref{eq:fond3}) and (\ref{eq:fond5})
with the claim of  \cite[Lemma 1]{Curtis2019},  we do not involve either the event $\nabla F(x_k)^Tg_k\ge 0$, nor
its probability.  
\vskip 5pt

The second result shows  the behaviour of the expected decrease of the optimality gap 
$F(x_{k+1})-F_*$ under the assumption that $F$ satisfies the PL condition, as well as the behaviour
of the expected  decrease of the average square 
norm of the gradient $\E \left[\frac{1}{K}\sum_{k=1}^K\|\nabla F(x_{k})\|^2 \right]$ in case $F$ is a 
general nonconvex  function.
Our analysis is carried out under the following assumption on the stochastic gradient.

\begin{ipotesi}\label{HMg}
There exists a positive constant $M_g$ such that, for all $k$,
$$
\Ek[\|\nabla F(x_k)-g_k\|^2] <M_g.
$$
\end{ipotesi}
This assumption was supposed to hold in the analysis of the second-order version of TRish carried out in \cite{Curtis2022}. 
{\color{black} It states that the variance of  $g_k$ is uniformly bounded and it can be nonzero at any stationary point of $F$. 
Such assumption includes the case where  $g_k$ is a stochastic approximation of $\nabla F(x_k)$ with additive noise, i.e., 
$g_k=\frac{1}{N}\sum_{i=1}^N (\nabla F_i(x_k)+\omega_i)$, each $\omega_i\in \R^{n}$ has mean zero and each component $(\omega_i)_j$ of $\omega_i$ has limited variance\footnote{{\color{black}
Trivially,  the triangle inequality implies $\|\nabla F(x_k)-g_k\|\le \sum_{i=1}^N \frac{1}{N}  \|\omega_i\|$. 
Then, by using  Jensen's inequality we obtain
$\|\nabla F(x_k)-g_k\|^2\le \left(\sum_{i=1}^N \frac{1}{N}  \|\omega_i\|\right)^2\le \frac{1}{N} \sum_{i=1}^N \|\omega_i\|^2$.
As a consequence, $\Ek[\|\nabla F(x_k)-g_k\|^2]\le \frac{1}{N} \sum_{i=1}^N \Ek[\|\omega_i\|^2]
=\frac{1}{N} \sum_{i=1}^N\sum_{j=1}^n \E_k[(w_i)_{j}^2]\le n\sigma$,
and Assumption \ref{HMg} holds.}}, i.e.,  $\Ek[\omega_i]=0$,  
$\Ek[(\omega_i)_j^2]\le \sigma$, for some $\sigma>0$, $1\le i\le N$, $1\le j\le n$.
In addition, Assumption \ref{HMg} holds when $g_k$ is un unbiased gradient estimate 
built by sampling, i.e., $g_k=\frac{1}{|S_k|}\sum_{i\in S_k}  \nabla F_i(x_k)$ with $S_k\subseteq\{1, \ldots, N\}$ and $|S_k|$ denoting the cardinality of the set $S_k$.
If the variance of each stochastic gradient $\nabla F_i(x_k)$ is equal and bounded by $M>0$, then the variance of  $g_k$ is bounded by $\frac{M}{|S_k|}$
(see \cite[\S 5.2]{nocedalsurvey}).
}
Assumption \ref{HMg} implies 
\begin{equation}\label{eq:ass_expgk}
\Ek[\|g_k\|^2]\leq M_1+M_2\|\nabla F(x_k)\|^2, \quad \forall \ k\geq 0,
\end{equation}
with $M_1=M_g$ and $M_2=1$.  Inequality \eqref{eq:ass_expgk} is generally assumed in the analysis of SG methods 
and also in the analysis of the first-order version of TRish presented in  \cite{Curtis2019}; clearly, it is a milder restriction
on the variance of $g_k$ than Assumption \ref{HMg}.
\vskip 5pt 

\begin{thm}\label{thm:gap}
Let Assumptions \ref{H1}, \ref{H3} and \ref{HMg} hold, $\beta$ as in (\ref{beta}), and suppose 
\begin{equation}\label{alpha}
\alpha< \frac{\gamma_2}{2\gamma_1^2 L} .
\end{equation}
\begin{enumerate}
\item 
If Assumption \ref{H2} holds and $\alpha< \min \left\{ \frac{\gamma_2}{2\gamma_1^2 L} , \frac{1 }{\mu \gamma_2}\right\}
$ then  
$$
\E[F(x_{k+1})]-F_* \overset{k\rightarrow\infty}{\longrightarrow}\frac{2\beta M_g}{\mu\gamma_2}.
$$
\item  If $f$ is a general function
$$
\E \left[\frac{1}{K}\sum_{k=1}^K\|\nabla F(x_{k})\|^2 \right] \overset{K\rightarrow\infty}{\longrightarrow} \frac{4\beta M_g}{\gamma_2}.
$$
\end{enumerate}
\end{thm}
{\bf Proof.}
1.  
By subtracting $F_*$ on both sides in (\ref{eq:fond5})  and using Assumption \ref{HMg}, we obtain
$$
\Ek[F(x_{k+1})]-F_*\le F(x_k)-F_*-\frac 1 2 \alpha \left(\gamma_2-\alpha \gamma_1^2 L \right )\|\nabla F(x_k)\|^2+
\alpha \beta M_g.
$$
Due to the bound on $\alpha$, we have $\gamma_2-\alpha\gamma_1^2 L>\gamma_2/2> 0$, and applying
the PL condition (\ref{eq:PL}) yields
\begin{equation}
\Ek[F(x_{k+1})]-F_*\le  \left(1-\frac 1 2 \alpha\mu \gamma_2 \right)(F(x_k)-F_*)+\alpha \beta M_g.
\end{equation}
By taking total expectations, letting $\xi=(1-\frac 1 2 \alpha\mu \gamma_2)$ and
proceeding by induction, we obtain
\begin{eqnarray*}
\E[F(x_{k+1})]-F_*&\le& \xi (\E[F(x_k)]-F_*)+\alpha \beta M_g\\
 &\leq &  \xi (\xi (\E[F(x_{k-1})]-F_*) +\alpha\beta M_g )+\alpha \beta  M_g  \\
& \leq & \cdots\leq   \xi^{k+1} (\E[F(x_0)]-F_*) + \sum_{j=0}^k\xi^k  \alpha \beta  M_g.
\end{eqnarray*}
The upper bound on  $\alpha$   and the relations $ \mu\le L$, $\gamma_2<\gamma_1$,  imply   $\xi \in (0, 1 )$ and  $\sum_{j=0}^\infty \xi^k=\frac{1}{1-\xi}$.
The claim follows by noting that $\frac{\alpha}{1-\xi}=\frac{2}{\gamma_2\mu}$.

2.
By applying Assumption \ref{HMg}, inequality (\ref{eq:fond5}) becomes
$$
\Ek[F(x_{k+1})]\le F(x_k)-\frac 1 4 \alpha \gamma_2\|\nabla F(x_k)\|^2+\alpha \beta M_g.
$$
By taking the total expectation, we get, for all $k$, 
$$
\E[\|\nabla F(x_k)\|^2]\le \frac{4 }{\alpha \gamma_2} (\E[F(x_{k})]-\E[F(x_{k+1})])+  \frac{4\beta M_g}{\gamma_2}.
$$
By summing this inequality for $k=1, \ldots, K$ and using that $F$ is bounded below by $F_*$, we obtain
$$
\E \left[\sum_{k=1}^K \|\nabla F(x_k)\|^2 \right]\le \frac{4}{\alpha\gamma_2} (F(x_{1})-F_*)+ K \frac{4\beta M_g}{\gamma_2} ,
$$
and the claim follows. $\Box$
\vskip 5pt 
Theorem \ref{thm:gap} holds for any value $\gamma_1/\gamma_2>1$ and differs in this respect from  
 \cite[Theorem 1, Theorem 5]{Curtis2019}, which prove, under  condition \eqref{eq:ass_expgk}, 
results on both the expected decrease of the 
optimality gap  and the expected decrease of the average square norm of the gradient.
In fact,  a restriction on the parameters $\gamma_1, \gamma_2$ of the form  
$\gamma_1-h_2(\gamma_1-\gamma_2)>0$, with $h_2\ge \frac 1 2 \sqrt{M_1}+\sqrt{M_2}$, is required
in \cite{Curtis2019}.
 Such restriction on the choice of $\gamma_1, \gamma_2$ is irrelevant
if $h_2\approx 1$, whereas it becomes stringent as $h_2$ increases. The need to assume a constraint on  $\gamma_1/\gamma_2$ in \cite{Curtis2019} depends on the more general assumption \eqref{eq:ass_expgk} imposed in place of Assumption \ref{HMg}.

We also observe that Theorem \ref{thm:gap} improves upon \cite[Theorem 4.1, Theorem 4.4]{Curtis2022}, which are also based on Assumption 
\ref{HMg},  in the values of 
the asymptotic expected 
optimality gap  and  average square norm of the gradient.
Specifically, the asymptotic expected optimality gap $\frac{2\beta M_g}{\mu \gamma_2}$
provided in Item 1 of Theorem \ref{thm:gap} 
is smaller than $\left(\frac{\gamma_1^2}{\gamma_2^2}-\frac{1}{2}\right)\frac{M_g}{\mu}$, due to the definition of $\beta$ in \eqref{beta} and the restriction of $\alpha$ in \eqref{alpha}, 
whereas the asymptotic expected optimality gap 
provided in \cite[Theorem 4.4]{Curtis2022} takes values  $4\left(\frac{\gamma_1^2}{\gamma_2^2}-\frac{1}{8}\right)\frac{M_g}{\mu}$.
As for the asymptotic  average square norm of the gradient, the value $\frac{4\beta M_g}{\gamma_2}$
provided in Item 2 of Theorem \ref{thm:gap} 
is smaller than $\left(2\frac{\gamma_1^2}{\gamma_2^2}-1\right)M_g$, due again to \eqref{beta} and \eqref{alpha}, 
whereas the asymptotic expected optimality gap 
provided in \cite[Theorem 4.1]{Curtis2022} takes values  $\left (8\frac{\gamma_1^2}{\gamma_2^2}-1\right)M_g$.

Finally, we note that the asymptotic expected optimality gap  and the average square norm of the gradient are proportional to $\beta$
in (\ref{beta}), which decreases as both $\gamma_1-\gamma_2$ and $\alpha$ decrease.
On the other hand,  the closer $\frac{\gamma_1}{\gamma_2}$ is to one, the closer the performance of
TRish algorithm becomes to that of the SG algorithm; likewise, the smaller $\alpha$, the slower the rate to achieve the limit value.

\vskip 5pt
We conclude our analysis by considering the case where condition 
(\ref{eq:ass_expgk}) holds with $M_1=0$ and for some positive $M_2$. Under this assumption, the optimality gap and the expected average square 
norm of the gradient eventually vanish.
This case covers the implementation of TRish considered in the next sections.

\vskip 5pt 
\begin{thm}\label{thm:nogap}
Let Assumptions \ref{H1} and \ref{H3} hold. Suppose that (\ref{eq:ass_expgk}) holds with $M_1=0$
 and that
$$
\Big (\frac{\gamma_2}{\gamma_1}\Big )^2 >1-\frac{1}{4 M_2}, \qquad 
\alpha<   \frac{1}{4\gamma_2LM_2} .
$$
\begin{enumerate}
\item 
If Assumption \ref{H2} holds and $\alpha< \min \left\{ \frac{1}{4\gamma_2L M_2}, \frac{2\gamma_2}{\mu \gamma_1^2}\right\}
$, then  
$$
\E[F(x_{k+1})]-F_* \overset{k\rightarrow\infty}{\longrightarrow} 0.
$$
\item  If $f$ is a general function, then
$$
\E \left[\frac{1}{K}\sum_{k=1}^K\|\nabla F(x_{k})\|^2 \right] \overset{K\rightarrow\infty}{\longrightarrow} 0.
$$
\end{enumerate}
\end{thm}
{\bf Proof.}
Recalling the definition (\ref{beta}) of  $\beta$, by \eqref{eq:fond3} we have 
\begin{eqnarray} 
\Ek[F(x_{k+1})]&\le& F(x_k)-\alpha \left (\frac{\gamma_1^2}{2\gamma_2} -\beta M_2\right) \|\nabla F(x_k)\|^2 \nonumber
\\
&\le& F(x_k)-\alpha \frac{\gamma_1^2}{4\gamma_2}   \|\nabla F(x_k)\|^2, \label{exp1}
\end{eqnarray}
where the last inequality follows from
$\beta<\frac{\gamma_1^2}{4\gamma_2 M_2}$ and $\frac{\gamma_1^2}{2\gamma_2} -\beta M_2>\frac{\gamma_1^2}{4\gamma_2}$. 

1. 
By subtracting $F_*$ on both sides  and using the PL condition (\ref{eq:PL}), we obtain
\begin{equation}
\Ek[F(x_{k+1})]-F_*\le  \left(1-2\alpha\mu \frac{\gamma_1^2}{4\gamma_2}  \right)(F(x_k)-F_*).
\end{equation}
By taking total expectations, letting $\xi=1- \alpha\mu \frac{\gamma_1^2}{2\gamma_2}$ and
proceeding by induction, we obtain
\begin{equation}
\E[F(x_{k+1})]-F_*\le \xi (\E[F(x_k)]-F_*)\le \ldots \le \xi^{k+1} (\E[F(x_0)]-F_*).
\end{equation}
Then, since $\xi\in (0,1)$ by definition, the claim follows.

2.
Taking total expectations on (\ref{exp1}) yields, for all $k$,
$$
\E[\|\nabla F(x_k)\|^2]\le \frac{4\gamma_2}{\alpha \gamma_1^2} (\E[F(x_{k})]-\E[F(x_{k+1})]).
$$
By summing this inequality for $k=1, \ldots, K$ and using that $F$ is bounded below by $F_*$, we obtain 
$$
\E \left[\sum_{k=1}^K \|\nabla F(x_k)\|^2 \right]\le \frac{4\gamma_2}{\alpha \gamma_1^2} (F(x_{1})-F_*),
$$
and the claim follows. $\Box$
\vskip 5pt
The claims of Theorem \ref{thm:nogap}
 are  similar to \cite[Theorem 1, Theorem 5]{Curtis2019} when one assumes $M_1=0$ in \eqref{eq:ass_expgk}, and they draw to conclusions that  are in accordance with the discussion in \cite{Curtis2019}.
In particular, the restriction on  $\frac{\gamma_1}{\gamma_2}$ in our theorem  can be written as 
$\gamma_1^2-4M_2(\gamma_1^2-\gamma_2^2)>0$; such a condition slightly affects the choice of $\gamma_1$ and $\gamma_2$ 
if $M_2$ is moderate, whereas it becomes more restrictive for large values of $M_2$, since it implies $\gamma_1\approx \gamma_2$ if $M_2\gg 0$. This issue is observed also in \cite[Theorem 1, Theorem 5]{Curtis2019}, where the requirement on $\gamma_1$ and $\gamma_2$ is $\gamma_1-\sqrt{M_2}(\gamma_1-\gamma_2)>0$. 

The requirements of Theorem \ref{eq:ass_expgk} and \cite[Theorem 1, Theorem 5]{Curtis2019} on the hyper-parameters $\alpha$, $\gamma_1$, $\gamma_2$, can be compared as follows. On the one hand, we observe that our upper bound on $\gamma_1/\gamma_2$ is more stringent than the one in \cite{Curtis2019}{\footnote{
The requirement on $\gamma_1$ and $\gamma_2$ in \cite[Theorem 1, Theorem 5]{Curtis2019}
is $\gamma_1-\sqrt{M_2}(\gamma_1-\gamma_2)>0$, i.e., 
$\frac{\gamma_2}{\gamma_1}< 1-\frac{1}{\sqrt{M_2}}$. 
It is easy to see that for $ M_2>1$ the condition $\Big (\frac{\gamma_2}{\gamma_1}\Big )^2 >1-\frac{1}{4 M_2}$
in Theorem \ref{thm:nogap}   is more restritive.
}.
On the other hand, our bounds on the steplength $\alpha$  improve upon those given in \cite[Theorem 1, Theorem 5]{Curtis2019} in the following respects.
Specifically, suppose $\gamma_2<\frac{1}{4}$ and  $M_2 >\max\{\gamma_1, \frac{\gamma_1^2\mu}{8\gamma_2^2 L}\}$. 
First, note that the condition $M_2 > \frac{\gamma_1^2\mu}{8\gamma_2^2 L}$ is not restrictive since $\frac{\gamma_2}{\gamma_1}<1$, $\mu\le L$.
Second, under Assumption \ref{H2}  the result in \cite[Eqn (17)]{Curtis2019} reads 
$\alpha\le \min\left\{\alpha^*, \alpha^\dagger\right\} := \min\left\{\frac{1}{\mu(\gamma_1-\sqrt{M_2}(\gamma_1-\gamma_2))},
\, \frac{\gamma_1-\sqrt{M_2}(\gamma_1-\gamma_2)}{\gamma_1LM_2}\right\}$ and 
the condition $\gamma_1<M_2$ easily yields $\alpha^\dagger =\min\left\{\alpha^*, \alpha^\dagger\right\}$.
Our convergence results are proved under the assumption $\alpha<\frac{1}{4\gamma_2L M_2}$ and 
 $\frac{1}{4\gamma_2L M_2}>\alpha^\dagger$
 is trivially satisfied whenever $\gamma_2<\frac{1}{4}$ since  $\frac{\gamma_1}{4\gamma_2}>\gamma_1$
while    $ (\gamma_1-\sqrt{M_2}(\gamma_1-\gamma_2) )< \gamma_1$.
Summarizing,  our upper bound on the steplength results larger than the one provided in \cite{Curtis2019}. For a general function, the discussion above still holds if $\gamma_2<\frac{1}{4}$.

We conclude this section by summarizing our results.
First, in Lemma \ref{lem:tec1} we provided results for the expected decrease $\Ek[F(x_{k+1})]-F(x_k)$ 
in a simplified way with respect to \cite{Curtis2019} and without invoking the probability of the
event $\nabla F(x_k)^T g_k\ge 0$. Second, in Theorem \ref{thm:gap} we provided asymptotic expected 
values under Assumption \ref{HMg} that are smaller than those  derived in \cite{Curtis2022}; 
{\color{black}we remark that Assumption \ref{HMg} was made in
\cite{Curtis2022} and  is stronger than the assumption made in \cite{Curtis2019} which is standard in the analysis of Stochastic Gradient methods}. 
Third,  in Theorem \ref{thm:nogap} we provided asymptotic expected values
under condition \eqref{eq:ass_expgk} with $M_1=0$, which is guaranteed by the implementation of TRish
considered in the following section; note that Theorem \ref{thm:nogap} improves the upper bound on the stepsize $\alpha$ under reasonable conditions on the parameters involved.
}

\section{TRish with adaptive sampling}\label{sec3}
In this section, we present a practical implementation of TRish in compliance with the theoretical assumptions of the previous section. Our implementation is based on the adaptive selection of the sample sizes 
for computing the stochastic estimators $g_k$  firstly proposed in \cite{Bollapragada2018}.
\subsection{Adaptive strategies for computing $g_k$}
Building the approximate gradient $g_k$ is a crucial and delicate issue in 
optimization methods based on random models.
In our case, due to the form (\ref{eq:problem}) of $F$, 
it is straightforward to consider a mini-batch gradient approximation. Thus, 
letting $S\subset\{1,\ldots, N\}$ be a sample, we denote  
\begin{equation}\label{eq:gradFS}
\nabla F_S(x)=\frac{1}{|S|}\sum_{i\in S}\nabla F_i(x),
\end{equation}
and we let  $g_k$ be of the form
\begin{equation}\label{eq:sampled_gradient}
g_k=\nabla F_{S_k}(x_k),
\end{equation}
which satisfies Assumption \ref{H3} if $S_k$ is uniformly chosen at random without replacement.

Forming $g_k$ leaves open the choice of the cardinality of the set $S_k$, which should be chosen in such a way that the assumptions made on $g_k$ in Section \ref{sec2} hold. 
A series of works \cite{BSV,BGMT, bcms,
Cartis-et-al, ChenMeniSche18,Bellavia-et-al-20, Bellavia-et-al-22, Byrd2012, Bollapragada2018, nocedalsurvey} 
present a dynamic choice of the sample size.  
In  \cite{Bellavia-et-al-20, Bellavia-et-al-22}, such  a choice is ruled by  
the combination of the Trust-Region and the Inexact Restoration
procedures. In \cite{BSV,BGMT, bcms,
Cartis-et-al, ChenMeniSche18}, the sample size is ruled by the 
Chebyshev inequality or the Bernstein concentration inequality.
In \cite{Byrd2012, Bollapragada2018}, the authors propose a  strategy for computing the sample size that enforces the theoretical properties specified in our analysis for $g_k$; in practice, they offer criteria for testing whether the current sample size should be maintained or increased. 
The common approach in \cite{Byrd2012, Bollapragada2018} is to estimate  the variance  in  the stochastic gradient
or in random variables depending on the stochastic gradient,
and increase the sample size if necessary. The heuristic relies on estimating
the variance by the sample variance;  if  such an estimation suggests an increase of the sample size,
it is assumed that the increase is small enough to consider reliable estimating the variance of the new sample by the variance of  the current one. Now we  summarize the strategies developed in  \cite{Byrd2012, Bollapragada2018}.

Condition (\ref{eq:start}) with $M_1=0$ can be enforced by adopting 
the strategy proposed by Bollapragrada et al.,  in  \cite{Bollapragada2018}, where 
$g_k$ is supposed to satisfy 
\begin{equation}\label{eq:inner_prod_test}
\Ek\left[(g_k^T\nabla F(x_k)-\|\nabla F(x_k)\|^2)^2\right]\leq \theta^2\|\nabla F(x_k)\|^4, \quad \text{for some }\theta>0,
\end{equation}
and  
\begin{equation}\label{eq:orthogonality_test}
\Ek\left[\left\|g_k-\frac{g_k^T\nabla F(x_k)}{\|\nabla F(x_k)\|^2}\nabla F(x_k)\right\|^2\right]\leq \nu^2\|\nabla F(x_k)\|^2, \quad \text{for some }\nu>0.
\end{equation}
Condition (\ref{eq:inner_prod_test}) represents an upper bound on the variance of $g_k^T\nabla F(x_k)$, while
condition (\ref{eq:orthogonality_test}) represents an upper bound on the variance of the component 
of $g_k$ that is orthogonal to the true gradient $\nabla F(x_k)$.
Under these conditions, (\ref{eq:ass_expgk}) holds with $M_1=0$ and $M_2=1+\theta^2+\nu^2$,
see \cite[Lemma 3.1]{Bollapragada2018}.
Since the above conditions are  too expensive to be implemented,  
the variances on the left-hand sides are approximated with the corresponding sample variances and the true gradient $\nabla F(x_k)$ with 
a sample gradient. Hence, the two conditions above are replaced as follows. Instead of condition (\ref{eq:inner_prod_test}), we use the approximate inner product test  
\begin{eqnarray}
& & {\rm Var}_{i\in S_k}(\nabla F_i(x_k)^Tg_k)  =\frac{1}{|S_k|-1}
\sum_{i\in S_k} \left( \nabla F_i(x_k)^Tg_k-\|g_k\|^2\right)^2, \nonumber \\ 
& & \frac{{\rm Var}_{i\in S_k}(\nabla F_i(x_k)^T g_k)}{ |S_k|}
\le \theta^2 \|g_k\|^4. \label{eq:inner_practical2}
\end{eqnarray}
Instead of condition (\ref{eq:orthogonality_test}), we use the  approximate orthogonality test
\begin{eqnarray}
& & {\rm Var}_{i\in S_k}\left (\nabla F_i(x_k)-\frac{\nabla F_i(x_k)^Tg_k}{\|g_k\|^2}g_k
\right) =\frac{1}{|S_k|-1} \sum_{i\in S_k}   \left\| \nabla F_i(x_k)-\frac{\nabla F_i(x_k)^Tg_k}{\|g_k\|^2}
\, g_k  \right\|^2 \nonumber \\ 
& & 
{\rm Var}_{i\in S_k}\left (\nabla F_i(x_k)-\frac{\nabla F_i(x_k)^Tg_k}{\|g_k\|^2}\, g_k
\right) \le \nu^2 \|g_k\|^2\label{eq:orthogonality_practical2}
\end{eqnarray}
If (\ref{eq:inner_practical2}) and (\ref{eq:orthogonality_practical2}) are satisfied, the sample size is 
kept unchanged, otherwise it is enlarged \cite[\S 4.3]{Bollapragada2018}.

\subsection{Algorithm implementation}
Our implementation is called TRish\_AS
(TRish with Adaptive Sampling) and combines TRish with the strategy based on conditions (\ref{eq:inner_practical2}) and (\ref{eq:orthogonality_practical2}). The resulting algorithm is in accordance with the requirements on $g_k$ specified in our theoretical analysis. Furthermore, it offers
an adaptive choice of $g_k$ and thus a rule for selecting the  parameters $\{|S_k|\}$.

We explain the main steps for computing the sample size in TRish\_AS. If conditions (\ref{eq:inner_practical2}) and (\ref{eq:orthogonality_practical2}) are satisfied, then  the sample size is 
kept unchanged. Otherwise, by using the simplifying assumptions in  \cite[\S 4.3]{Bollapragada2018}, a sample size larger than 
$S_k$ is computed as 
\begin{equation}\label{eq:sk_ok}
\min\left \{ \max \left\{\ \left\lceil \frac{{\rm Var}_{i\in S_k}(\nabla F_i(x_k)^T g_k)}{
\theta^2\|g_k\|^4}\right\rceil, \left\lceil 
\frac{{\rm Var}_{i\in S_k}\left (\nabla F_i(x_k)-\frac{\nabla F_i(x_k)^T g_k}{\|g_k\|^2}\, g_k
\right) }{\nu^2 \|g_k\|^2}\right\rceil
 \right\}, N \right \}.
\end{equation}
Finally, we  use the sample control in the noisy regime fully described in \cite[\S 4.2]{Bollapragada2018}
and sketched at Step 5 of Algorithm \ref{algo:trish_practical}. 
If the sample sizes do not change along $r$ iterations, the 
average of the most $r$ recent sample gradients
$g_{\textrm{avg}}=\frac 1 r \sum_{j=k-r+1}^k\nabla F_{S_j}(x_j)$  is computed.
Then, if the norm of $g_{\textrm{avg}}$ is small compared to the norm of $g_k$,
(\ref{eq:inner_practical2})  and (\ref{eq:orthogonality_practical2}) are tested 
 using $g_{\textrm{avg}}$ instead of $\nabla F_{S_k}(x_k)$; in case of failure of one of such conditions, the sample size is  increased  using (\ref{eq:sk_ok}) and $g_{\textrm{avg}}$ instead of $g_k$. We report TRish\_AS in the following Algorithm \ref{algo:trish_practical}.

\begin{algorithm}[h!]
\caption{{\bf TRish\_AS: practical implementation}}\label{algo:trish_practical}
\vspace{2mm}
\noindent Choose an initial iterate $x_0\in\mathbb{R}^n$, a steplength $\alpha>0$, a sample size $S_0\in \{1, \ldots, N\}$ 
an integer $r>0$ and real parameters $0<\gamma_2<\gamma_1$, $\theta>0$, $\nu>0$.\\
\vskip 1pt
Set $k=0$. Compute $g_0=\nabla F_{S_0}(x_0)$.
\vskip 2pt\noindent
Repeat until a convergence test is satisfied
\begin{enumerate}
\item Compute
\begin{equation}\label{eq:trish_stepB}
x_{k+1}=x_k-\begin{cases}
\gamma_1 \alpha g_k, \quad &\text{if }\|g_k\|\in[0,\frac{1}{\gamma_1})\\
\alpha \displaystyle \frac{g_k}{\|g_k\|}, \quad &\text{if }\|g_k\|\in[\frac{1}{\gamma_1},\frac{1}{\gamma_2}] \\
\gamma_2 \alpha g_k, \quad &\text{if }\|g_k\|\in(\frac{1}{\gamma_2},+\infty[.
\end{cases}
\end{equation}
\item Set $k=k+1$.
\item Set  $|S_k|=|S_{k-1}|$, choose a new sample $S_k$  and form $g_k=\nabla F_{S_k}(x_k)$.
\item  If (\ref{eq:inner_practical2}) or (\ref{eq:orthogonality_practical2}) is not satisfied then
\\
\hspace*{10pt} update $|S_k|$ using (\ref{eq:sk_ok}), choose a new sample $S_k$ and form
$g_k=\nabla F_{S_k}(x_k)$. 
\item If $|S_k|=|S_{k-1}|=\cdots =|S_{k-r}|$ then
\\
\hspace*{10pt} compute $g_{\textrm{avg}}=\frac 1 r \sum_{j=k-r+1}^k\nabla F_{S_j}(x_j)$.
\\
\hspace*{10pt}
If $\|g_{\textrm{avg}}\|<\gamma\|\nabla F_{S_k}(x_k)\| $ then
\\
\hspace*{25pt} If (\ref{eq:inner_practical2}) or (\ref{eq:orthogonality_practical2}) is not satisfied 
 using $g_{\textrm{avg}}$ instead of $\nabla F_{S_k}(x_k)$ then 
\\
\hspace*{40pt} update $|S_k|$ using (\ref{eq:sk_ok}) and $g_{\textrm{avg}}$ instead of $g_k$;
\\
\hspace*{40pt} choose a new sample $S_k$  and form $g_k=\nabla F_{S_k}(x_k)$.
\end{enumerate}
\end{algorithm}

\section{Numerical experience}\label{sec4}
We now present the numerical performance of Algorithm \ref{algo:trish_practical}. Our goal is to show that
the overall performance of TRish  may benefit from this adaptive approach, instead of using a prefixed and constant sample size as done in \cite{Curtis2019, Curtis2022}.

\subsection{Logistic Regression}\label{sec:regression}
We consider datasets with feature vectors $z_i\in \R^n$, labels $y_i\in \{-1, 1\}$, $i=1, \ldots, N$,
and address the problem of binary classification by minimizing the logistic regression function
$$
F(x)=\frac 1 N \sum_{i=1}^N \log \left(1+e^{-y_i(x^Tz_i)}\right).
$$
The datasets considered are available in the LIBSVM repository \cite{libsvm}.
A testing dataset $\{ (\bar z_i, \bar y_i) \}_{i=1}^{\bar{N}}$ is used for evaluating the accuracy in classifying
unseen data.

Our experiments were carried out as in \cite{Curtis2019}. 
We considered $60$ parameter settings $(\alpha, \gamma_1, \gamma_2)$ defined as follows.
The stepsize takes value
$\alpha\in\{0.1, 10^{-\frac 1 2},1, 10^{\frac 1 2}, 10 \}$. The parameters $\gamma_1, \gamma_2$ take value
$\gamma_1\in\{\frac 4 G, \frac 8 G , \frac{16}{G}, \frac{32}{G}\}$, $\gamma_2\in\{\frac{1}{2G}, \frac{1}{G}, \frac{2}{G}\}$
with $G$ being the average norm of stochastic gradients generated by SG run for one epoch 
with stepsize $\alpha=0.1$ and sample size equal to 64. 
A complete set of experiments amounts to applying both TRish and TRish\_AS with 
the 60  different combinations  $(\alpha, \gamma_1, \gamma_2)$. Both algorithms were run for one epoch. TRish was run using a sample size equal to 64, see \cite{Curtis2019}, while the initial sample size for TRish\_AS was $S_0=\min\left\{32, \left\lceil \frac{N}{100}\right\rceil\right\}$. Regarding the implementation of condition \eqref{eq:sk_ok} in TRish\_AS, we fixed  $\theta=0.9$ and $\nu=5.84$, due to their statistical significance (see \cite[\S 4.4]{Bollapragada2018}) and  $r=10$. In a few cases, the evaluation of  (\ref{eq:sk_ok}) in finite precision was harmful, due to very large or small norm of $g_k$, and affected by underflow or overflow; in such occurrences the sample size was not updated.

In the remainder of the section, we report the following results. 
For each  parameter setting  $(\alpha, \gamma_1, \gamma_2)$ 
we measure, averaging over 50 runs, the testing accuracy defined as the fraction of the testing set correctly classified
and, limited to TRish\_AS, the  sample size reached at termination.
Then,  for each test  problem and parameter setting we plot the average testing accuracy 
reached by TRish (star symbol) and TRish\_AS (circle symbol). For the sake of clarity, we also provide the same plots 
separately for each fixed value of $\alpha$ and all tested values of $\gamma_1, \gamma_2$, i.e., we 
focus  on the  12 parameter  settings sharing the same value of $\alpha$.
For each test problem we also present  the plots on the average training loss and testing accuracy versus  the average number
of effective gradient evaluations (EGE), that is the number of full gradient evaluations. Finally,  in a table we summarize the statistics 
corresponding to the best run for TRish (first row), i.e., the case   where  
TRish reached the best average accuracy, and  the statistics 
corresponding to the best run for TRish\_AS (second row).
On each row, we report  the  parameter setting $(\alpha, \gamma_1, \gamma_2)$  that provided the
best run and, for such a triplet, the average testing accuracy obtained by the two methods and the average  final  
sample size used by TRish\_AS; thus, it is possible to compare the results obtained by the two algorithms  
with the same parameter setting. Correspondingly to such runs we also give 
statistics on the steps taken by TRish and TRish\_AS
among the normalized step and the SG-type step (Cases 1--3).

\subsubsection{{\bf a1a}}
\begin{figure}[t!]
\begin{center}
\begin{tabular}{c}
\includegraphics[scale = 0.4]{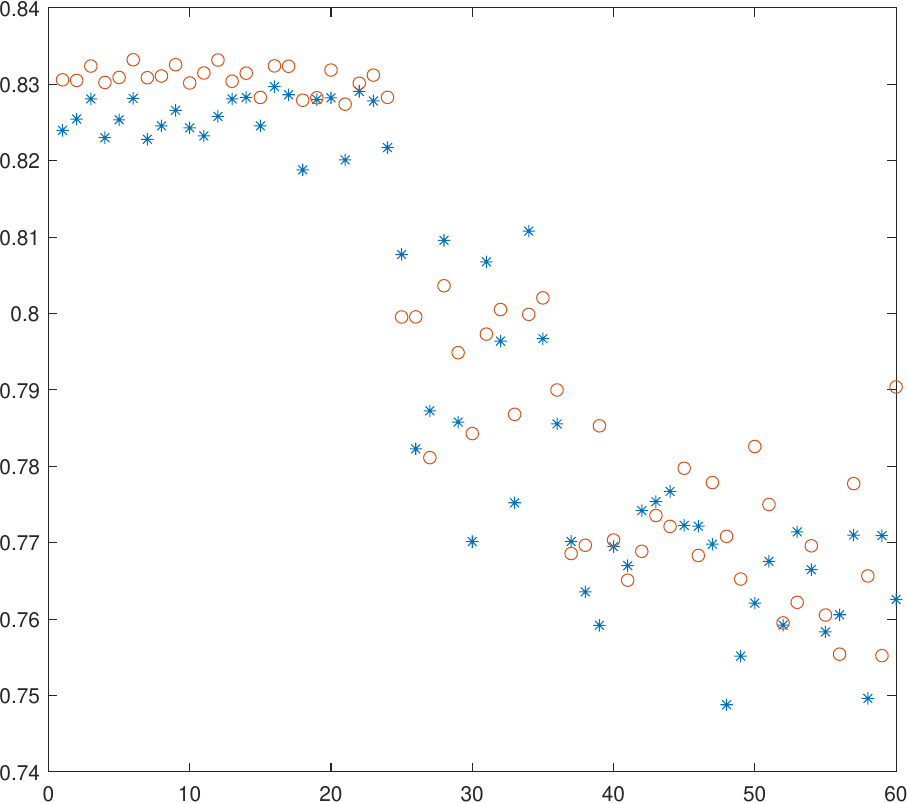}
\end{tabular}
\caption{Dataset  a1a. Average testing accuracy at termination. 
TRish: symbol  ``*'', TRish\_AS: symbol ``o''. }
\label{fig1:a1a}
\end{center}
\end{figure}

\begin{figure}[t!]
\begin{center}
\begin{tabular}{c}
\includegraphics[scale = 0.4]{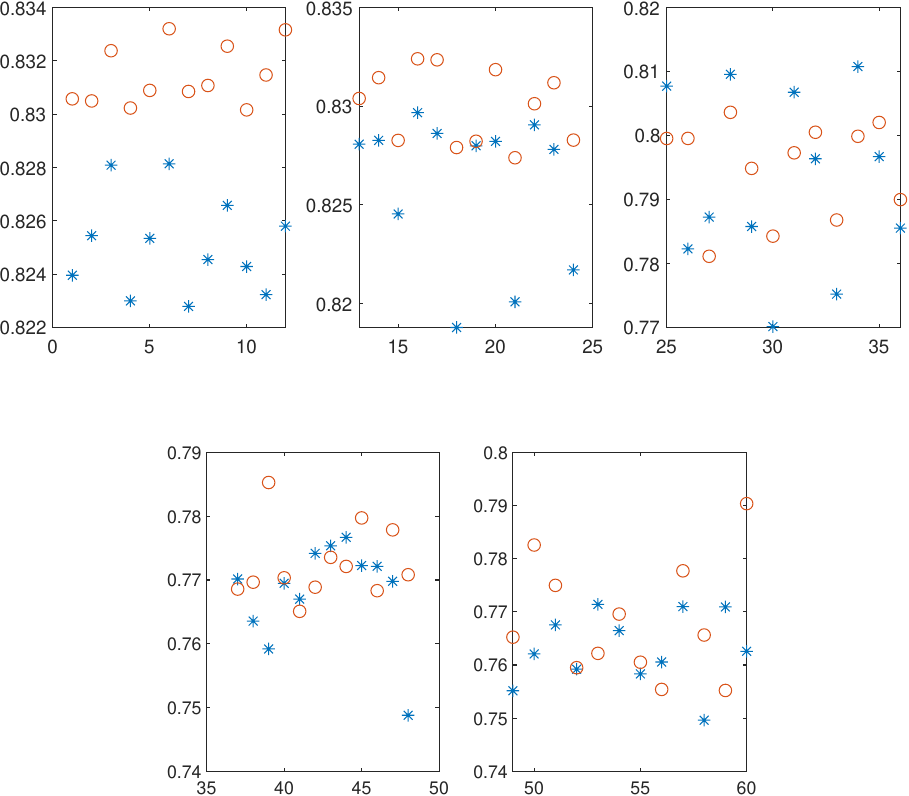}
\end{tabular}
\caption{Dataset a1a. Average testing accuracy at termination. 
Runs with fixed steplength $\alpha$ and varying $\gamma_1, \gamma_2$.
Top, from left to right:  $\alpha=10^{-1}$, $\alpha=10^{-\frac 1 2}$, $\alpha =1$. Bottom, from left to right: $\alpha=10^{\frac 1 2}$, $\alpha=10$. 
TRish: symbol  ``*'', TRish\_AS: symbol ``o''. }
\label{fig1:a1aplot4}
\end{center}
\end{figure}

We used a training set consisting of $N=1605$ points, features vectors have dimension $n=123$.
The testing set consists of $\bar N=29351$ points.
We performed the runs using the value $G=0.3477$ for the setting of $\gamma_1$ and $\gamma_2$. 
In Figure \ref{fig1:a1a}  we plot the average testing accuracy reached by TRish
and by TRish\_AS, while in  Figure \ref{fig1:a1aplot4} each plot displays the average testing accuracies 
reached for specific value of  $\alpha$ and all tested values of $\gamma_1, \gamma_2$. 

\begin{figure}[t!]
\begin{center}
\begin{tabular}{c}
\includegraphics[scale = 0.4]{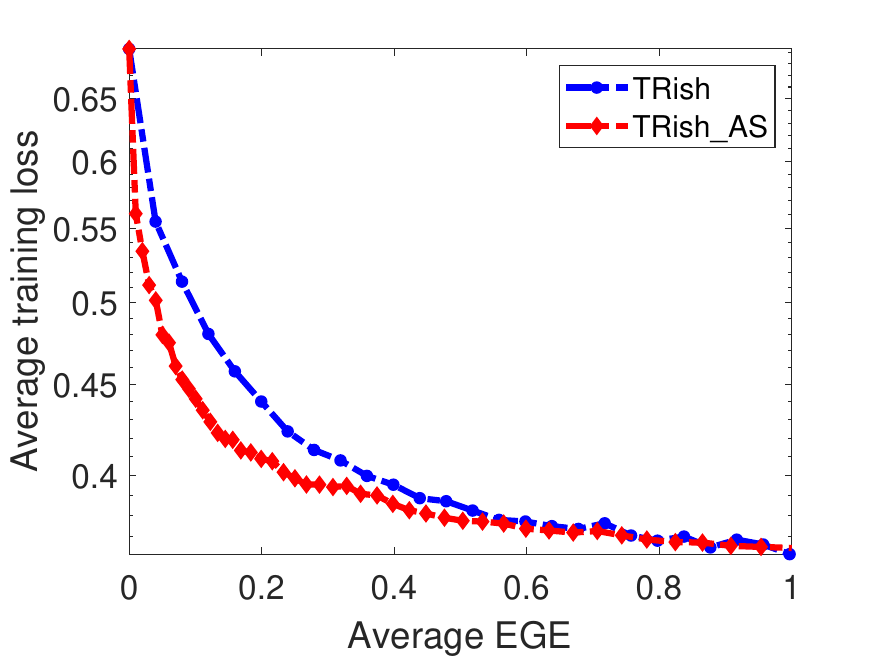}
\includegraphics[scale = 0.4]{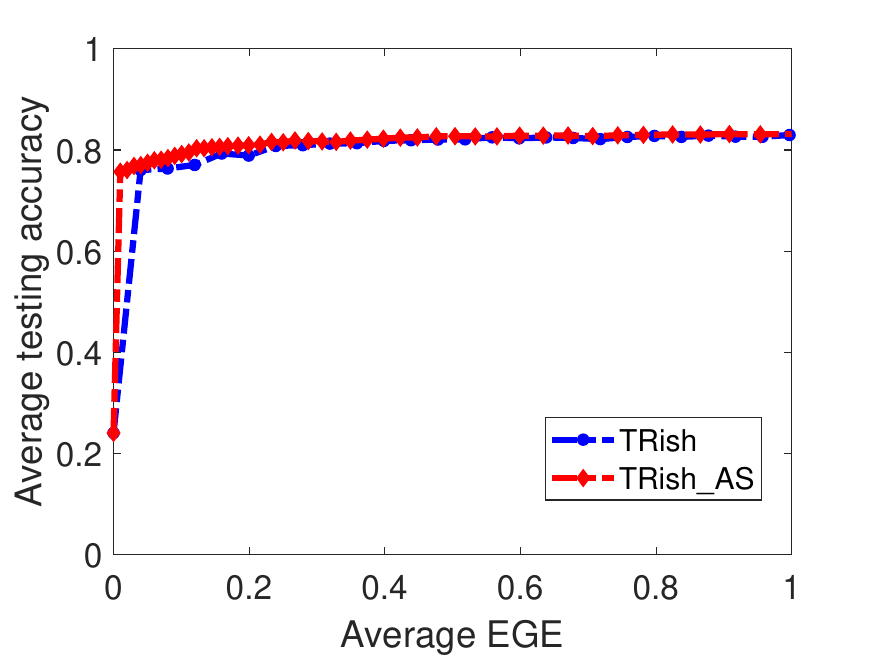}
\end{tabular}
\caption{Dataset a1a. Average training loss (left) and testing accuracy (right) versus one epoch for 
TRish (blue) and  TRish\_AS (red).}
\label{fig4:a1a}
\end{center}
\end{figure}

The results in Figure \ref{fig1:a1a} show that the accuracy attained with TRish\_AS is higher than the accuracy 
obtained by  TRish  in 46 cases out of 60. The minimum and maximum averaged sample sizes at termination of TRish\_AS 
were 
60 and 186 respectively.  We also observe in  Figure \ref{fig1:a1aplot4}  that when  the
steplengths  $\alpha=10^{-1}, 10^{-\frac 1 2 }$  are used, 
TRish\_AS is quite insensitive to the choice of $\gamma_1$ and $\gamma_2$  and TRish\_AS  is more accurate 
than TRish in classification. 
In Figures \ref{fig4:a1a}   we show the behaviour of the average training loss and testing accuracy
versus  the average number
of effective gradient evaluations (EGE) for the two methods equipped with their own best parameters settings.   

Table \ref{tab:a1a} concerns the best (averaged) run and  
shows that TRish\_AS outperforms TRish for its best parameters setting, as well as for the one which provides the best averaged  run of TRish. Moreover, the  average final sample sizes of TRish\_AS
are considerably  small,  as they represent the 7.7\% of the whole sample for the first parameter combination, and the 4.7\% for the second parameter combination, respectively.
We also give statistics regarding the  step taken along the iterations among  the normalized step (Case 2) and the SG-type steps (Cases 1 and 3).
Focusing on the runs with parameters given in  the second row of Table \ref{tab:a1a}, on average TRish selected the normalized step (Case 2) $33\%$   times
and the SG-type step (Case 3) $67\%$  times, while on average TRish\_AS selected the normalized step (Case 2) $19\%$ times
and the SG-type step (Case 3) $81\%$   times.

\begin{table}[h!]
\small
\centering
\begin{tabular}{ l|c|c|c|c}
\multicolumn{ 1}{c|}{} &\multicolumn{ 1}{c|}{$(\alpha, \gamma_1, \gamma_2)$}  & \multicolumn{ 2}{c |}{Testing accuracy} & {Average final $|S_k|$}\\ \hline
  &  & \multicolumn{ 1}{c|}{  TRish} & \multicolumn{ 1}{c|}{  TRish\_AS}  & { TRish\_AS }  \\ \hline
 Best run of TRish & $(  0.3162 ,  23.9593  ,  1.4975  ) $ &  0.8297   & 0.8324 & 124  \\
Best run of TRish\_AS  & $(0.1000,   23.9593 ,   5.9898)$ &   0.8281  &  0.8332 &  75   \\
\end{tabular}  
\caption{Dataset a1a.  First row:   parameter setting corresponding to the best run for TRish in terms of testing accuracy, average testing accuracy of TRish and TRish\_AS at termination with the triplet,  average sample size of TRish\_AS at termination  with the triplet.
 Second row:   parameter setting corresponding to the best run for TRish\_AS in terms of testing accuracy, average testing accuracy of TRish and TRish\_AS at termination with the triplet,  average sample size of TRish\_AS at termination  with the triplet.   }\label{tab:a1a}
\end{table} 

\subsubsection{{\bf w1a}}
This  training data set  consists of $N=2477$ points, features vectors have dimension $n=300$.
The testing set consists of $\bar N=47272$ points.
Following \cite{Curtis2019} we performed the runs using the value $G=0.0887$ for the setting of $\gamma_1$ and $\gamma_2$. 

\begin{figure}[t!]
\begin{center}
\begin{tabular}{c}
\includegraphics[scale = 0.4]{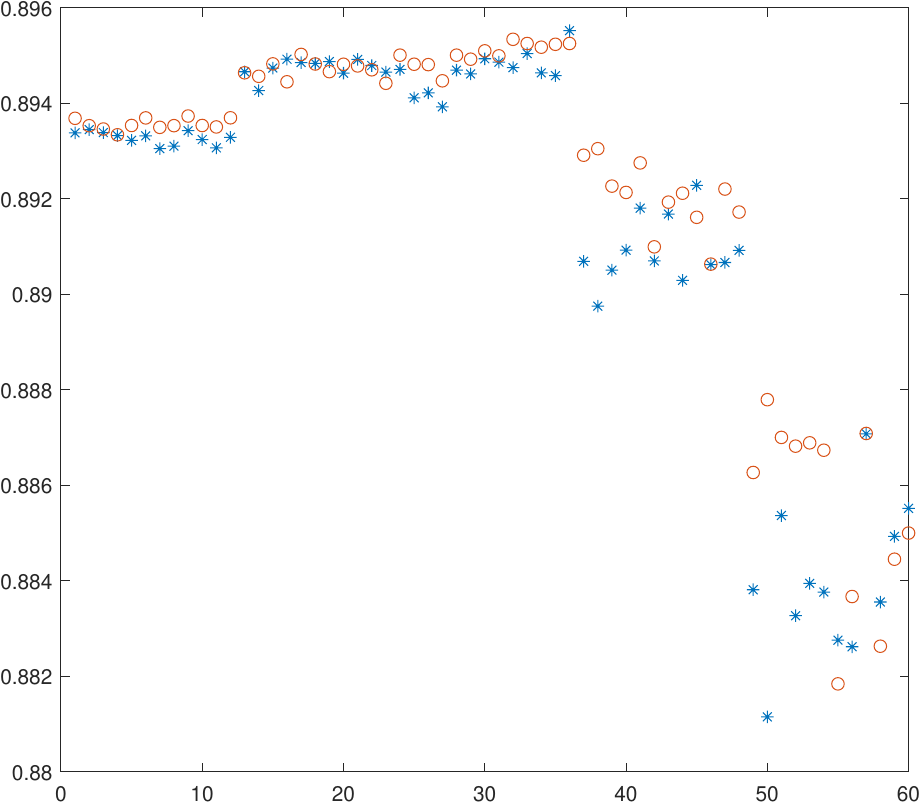}
\end{tabular}
\caption{Dataset  w1a. Average testing accuracy at termination. 
TRish: symbol  ``*'', TRish\_AS: symbol ``o''. }
\label{fig1:w1a}
\end{center}
\end{figure}

\begin{figure}[t!]
\begin{center}
\begin{tabular}{c}
\includegraphics[scale = 0.4]{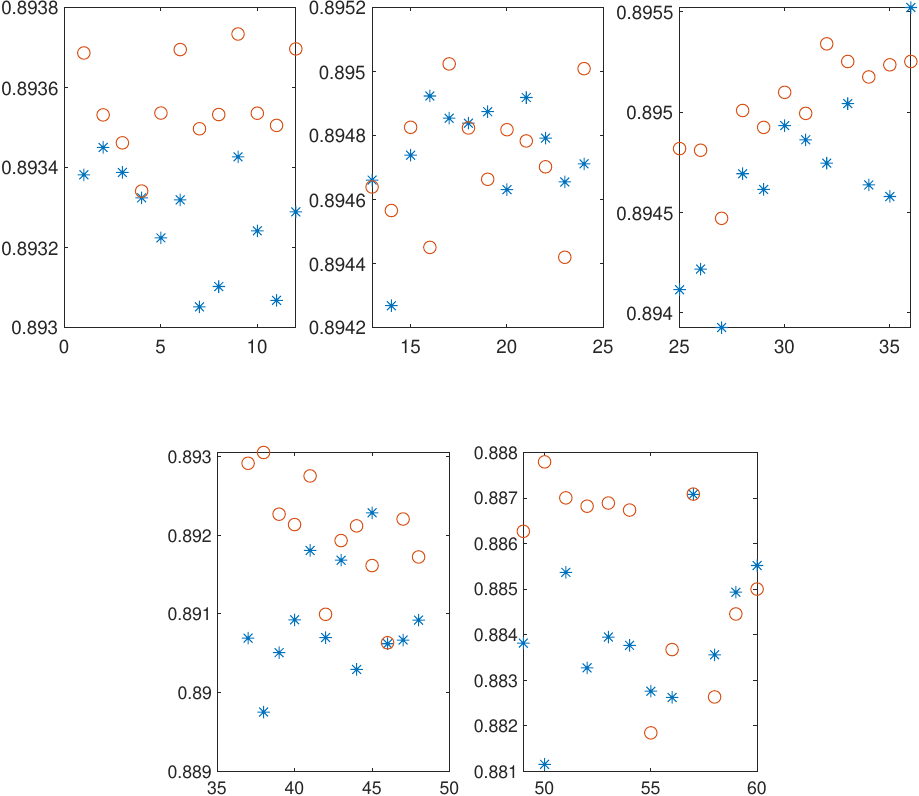}
\end{tabular}
\caption{Dataset  w1a. Average testing accuracy at termination. 
Runs with fixed steplength $\alpha$ and varying $\gamma_1, \gamma_2$.
Top, from left to right:  $\alpha=10^{-1}$, $\alpha=10^{-\frac 1 2}$, $\alpha =1$. Bottom, from left to right: $\alpha=10^{\frac 1 2}$, $\alpha=10$. 
TRish: symbol  ``*'', TRish\_AS: symbol ``o''. }
\label{fig1_w1aplot4}
\end{center}
\end{figure}

\begin{figure}[t!]
\begin{center}
\begin{tabular}{c}
\includegraphics[scale = 0.4]{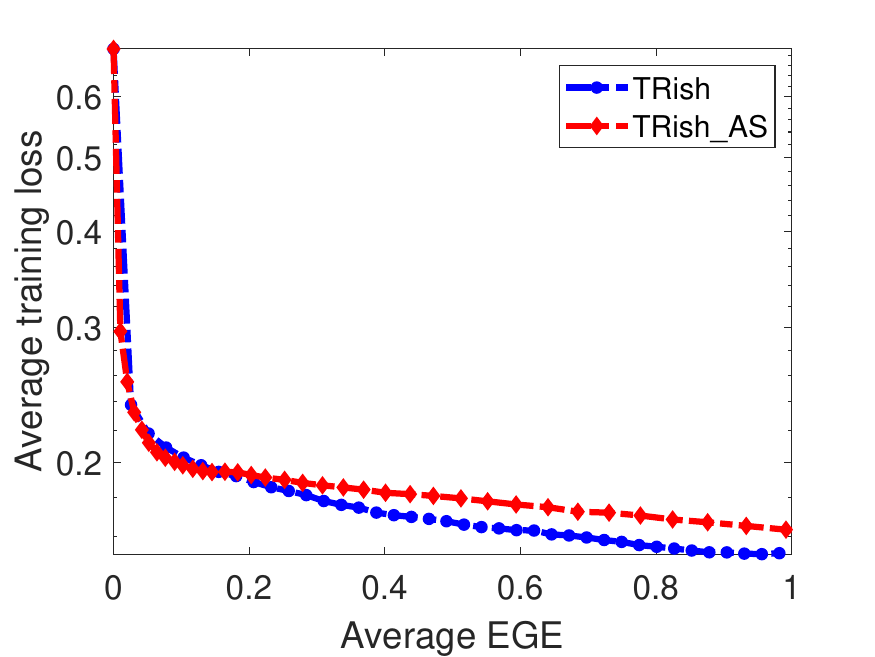}
\includegraphics[scale = 0.4]{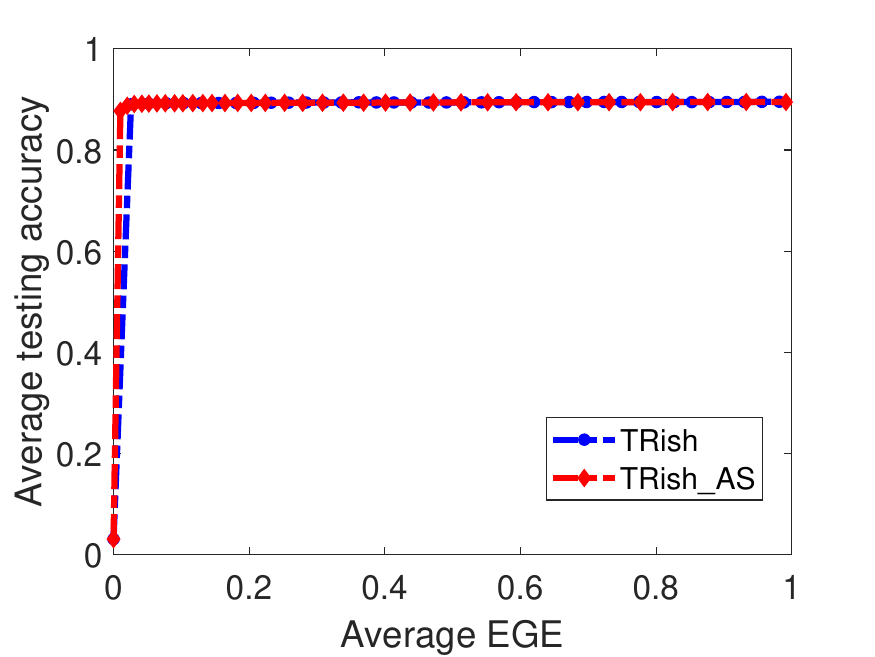}
\end{tabular}
\caption{Dataset w1a. Average training loss (left) and testing accuracy (right) versus one epoch for 
TRish (blue) and  TRish\_AS (red). }
\label{fig4:w1a}
\end{center}
\end{figure}

Figure \ref{fig1:w1a} shows that TRish\_AS  obtained the highest accuracy 
in 47 cases out of 60. The minimum and maximum average sample sizes of TRish\_AS at termination were 
67 and 343 respectively.   
Figure \ref{fig1_w1aplot4} clarifies that the accuracies achieved by the two methods are similar, and that 
the performance of TRish\_AS is only slightly affected by the choice of parameters $\gamma_1$ and $\gamma_2$.
In Figures \ref{fig4:w1a}   we show the behaviour of the average training loss and testing accuracy
versus the average number of  EGE for the methods equipped with their own best parameters settings. 
We observe that  TRish reaches smaller values of the training loss, while the two methods achieves almost the same accuracy.
 
The good behaviour of TRish\_AS is confirmed in Table 
\ref{tab:w1a} as  TRish\_AS is comparable to   TRish in  the best performance of the latter. 
In the reported runs the sample size at termination is  around the 7.4\%  and 9.8\% of $N$ respectively. When the parameters were selected as in the second row of Table \ref{tab:w1a}, TRish computed the normalized step (Case 2)
$85\%$  times and the SG-type step (Case 1 and 3) $15\%$   times on average, whereas TRish\_AS computed 
the normalized step (Case 2)
$83\%$ times and the SG-type step (Case 1 and 3) $17\%$  times on average.
 
\begin{table}[h!]
\small
\centering 
\begin{tabular}{ l|c|c|c|c}
\multicolumn{ 1}{c|}{} &\multicolumn{ 1}{c|}{$(\alpha, \gamma_1, \gamma_2)$}  & \multicolumn{ 2}{c |}{Testing accuracy} & {Average final $|S_k|$}\\ \hline
  &  & \multicolumn{ 1}{c|}{  TRish} & \multicolumn{ 1}{c|}{  TRish\_AS}  & { TRish\_AS }  \\ \hline
 Best run of TRish   & $(1.0000,  105.8551  ,  6.6159)$      &0.8955 &   0.8953  & 183\\
Best run of TRish\_AS  &   $(1.0000,   52.9276,    3.3080)$ &     0.8947 &   0.8953  &242 \\
\end{tabular}  
\caption{Dataset  w1a.   First row:   parameter setting corresponding to the best run for TRish in terms of testing accuracy, average testing accuracy of TRish and TRish\_AS at termination with the triplet,  average sample size of TRish\_AS at termination  with the triplet.
 Second row:   parameter setting corresponding to the best run for TRish\_AS in terms of testing accuracy, average testing accuracy of TRish and TRish\_AS at termination with the triplet,  average sample size of TRish\_AS at termination  with the triplet.}\label{tab:w1a}
\end{table} 

\subsubsection{{\bf rcv1}}
This  training data set  consists of $N=20242$ points, features vectors have dimension $n=47236$.
The testing set consists of $\bar{N}=677399$ points.
Following \cite{Curtis2019} we performed the runs using the value $G=0.0497$ for the setting of $\gamma_1$ and $\gamma_2$. 
{Figure \ref{fig1:rcv1} shows that TRish\_AS outperforms TRish for almost all parameters combinations; more specifically, TRIsh\_AS is more accurate than TRish $59$ out of $60$ runs. However, both methods perform quite well and the lowest fraction of correctly classified features is  approximately $0.9303$ attained by TRish with  $\alpha=10^{-1}$, see Figure \ref{fig1_rcv1plot4}. Figure \ref{fig4:rcv1} shows the decrease of the average training loss and testing accuracy
versus the average number
of effective gradient evaluations  for the methods equipped with their own best parameters settings. The best run of each procedure correctly classified  95\% of instances in the testing set as reported in Table \ref{tab:rcv1}.
The  minimum averaged sample size at termination of TRish\_AS was $32$, while the maximum was $178$ corresponding to  0.87\% of $N$.
Finally, considering the parameters setting at the second row of Table \ref{tab:rcv1}, we note that TRish selected the normalized step (Case 2) $75\%$   times 
and the SG-type step (Case 1 and 3) $25\%$  times on average, whereas TRish\_AS selected the normalized step (Case 2) $69\%$ times
and the SG-type step (Case 3) $31\%$   times on average.

\begin{figure}[t!]
\begin{center}
\begin{tabular}{c}
\includegraphics[scale = 0.4]{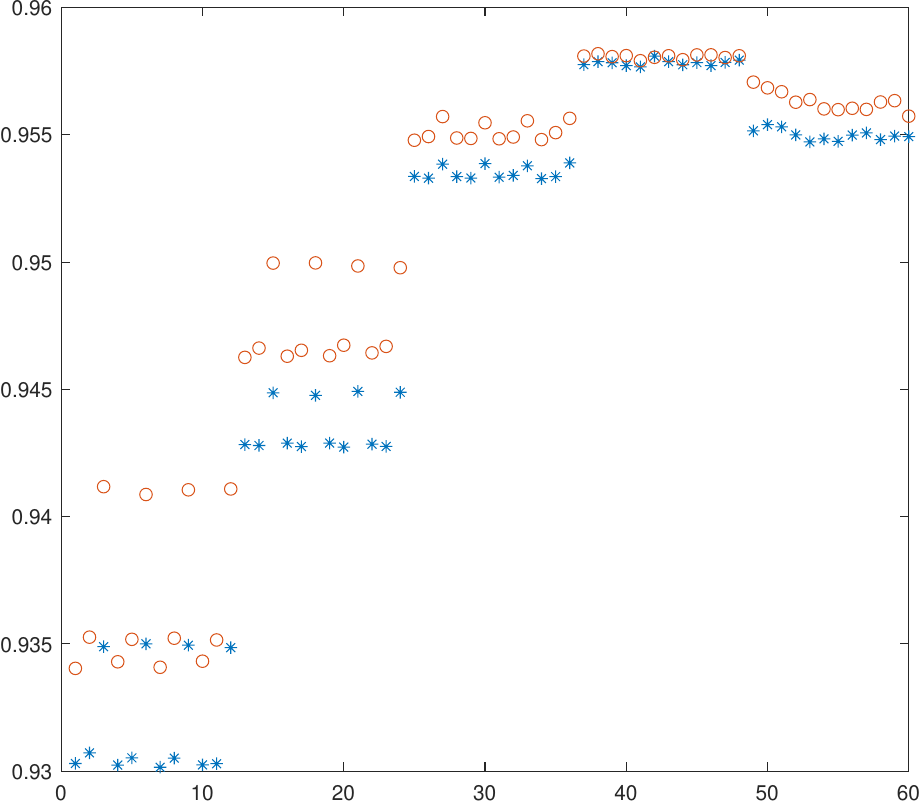}
\end{tabular}
\caption{Dataset   rcv1. Average testing accuracy at termination. 
TRish: symbol  ``*'', TRish\_AS: symbol ``o''. }
\label{fig1:rcv1}
\end{center}
\end{figure}

\begin{figure}[t!]
\begin{center}
\begin{tabular}{c}
\includegraphics[scale = 0.4]{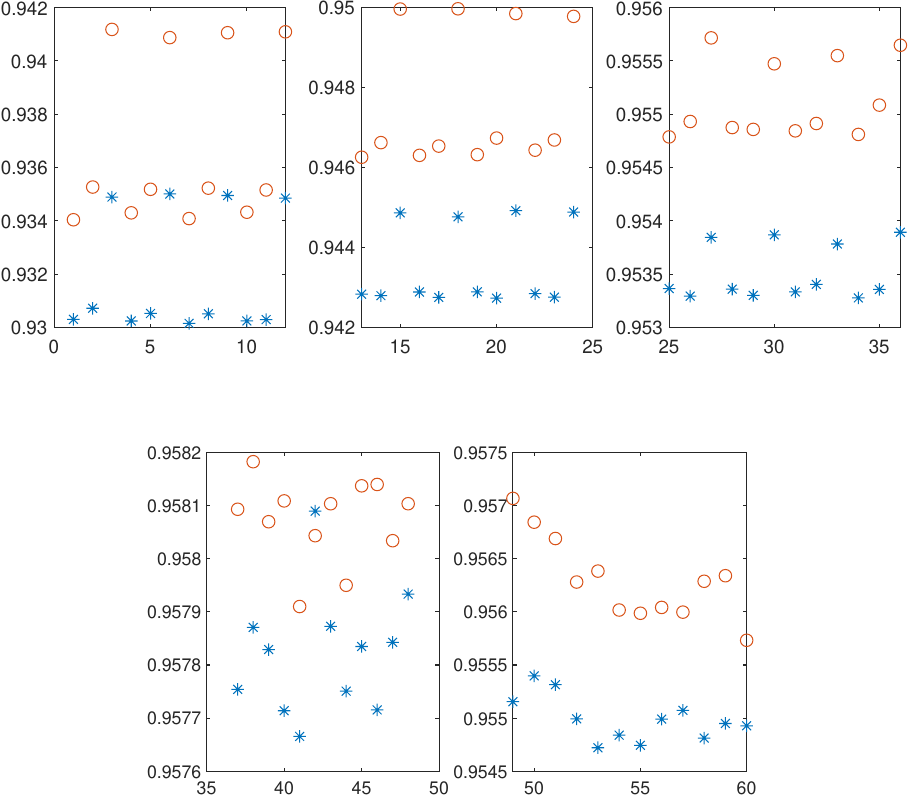}
\end{tabular}
\caption{Dataset rcv1. Average testing accuracy at termination. 
Runs with fixed steplength $\alpha$ and varying $\gamma_1, \gamma_2$. Top, from left to right:  $\alpha=10^{-1}$, $\alpha=10^{-\frac{1}{2}}$, $\alpha =1$. Bottom, from left to right: $\alpha=10^{\frac{1}{2}}$, $\alpha=10$.
TRish: symbol  ``*'', TRish\_AS: symbol ``o''. }
\label{fig1_rcv1plot4}
\end{center}
\end{figure}

\begin{figure}[t!]
\begin{center}
\begin{tabular}{c}
\includegraphics[scale = 0.4]{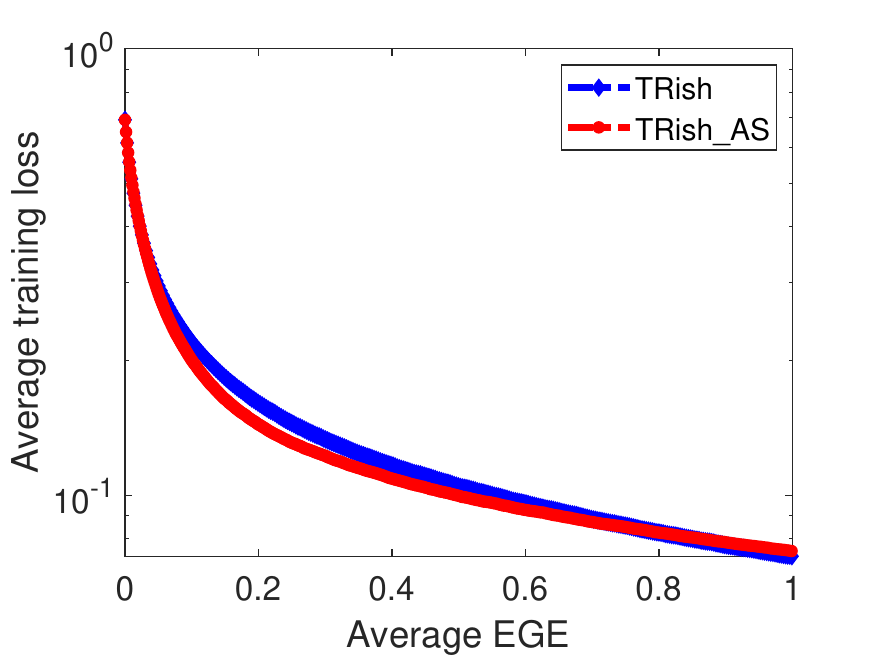}
\includegraphics[scale = 0.4]{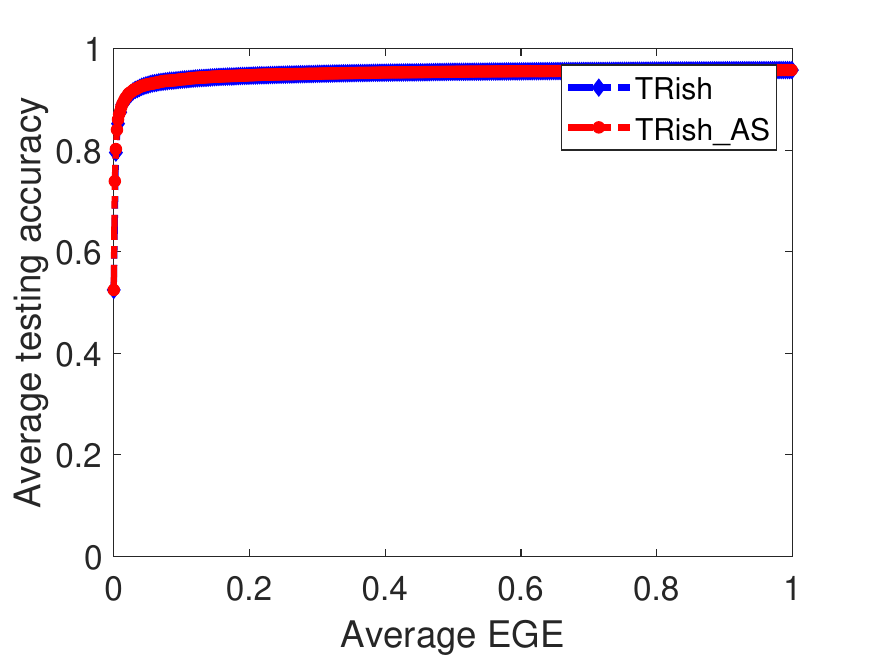}
\end{tabular}
\caption{ Dataset rcv1. Average training loss (left) and testing accuracy (right) versus one epoch for 
TRish (blue) and  TRish\_AS (red). }
\label{fig4:rcv1}
\end{center}
\end{figure}

\begin{table}[h!]
\small
\centering
\begin{tabular}{ l|c|c|c|c}
\multicolumn{ 1}{c|}{} &\multicolumn{ 1}{c|}{$(\alpha, \gamma_1, \gamma_2)$}  & \multicolumn{ 2}{c |}{Testing accuracy} & {Average final $|S_k|$}\\ \hline
  &  & \multicolumn{ 1}{c|}{  TRish} & \multicolumn{ 1}{c|}{  TRish\_AS}  & { TRish\_AS }  \\ \hline
 Best run of TRish   & $(3.1623,  121.3961,   30.3490)$  &  0.9581  &  0.9580 & 116\\
Best run of TRish\_AS  &   $(3.1623 ,  60.6980,   15.1745 )$ &   0.9579 &   0.9582   &  109  \\
\end{tabular}  
\caption{Dataset  rcv1.   First row:   parameter setting corresponding to the best run for TRish in terms of testing accuracy, average testing accuracy of TRish and TRish\_AS at termination with the triplet,  average sample size of TRish\_AS at termination  with the triplet.
 Second row:   parameter setting corresponding to the best run for TRish\_AS in terms of testing accuracy, average testing accuracy of TRish and TRish\_AS at termination with the triplet,  average sample size of TRish\_AS at termination  with the triplet.} \label{tab:rcv1}
\end{table} 

\subsubsection{{\bf Cina}}
This training data set  consists of $N=10000$ points, features vectors have dimension $n=132$.
The testing set consists of $\bar N=6033$ points.  The dataset was normalized such that each
feature vector has components in the interval $[0,1]$.
We performed the runs using the value $G=0.2517$ for the setting of $\gamma_1$ and $\gamma_2$. 

TRish\_AS   achieves the highest accuracy in classification
in 17 cases out of 60, see Figures  \ref{fig1:cina} and  \ref{fig1_cinaplot4}. 
The accuracy is very similar for the two smallest steps, $\alpha=10^{-1}$ and $\alpha=10^{-\frac 1 2 }$, 
which produce the highest accuracy. Overall the largest gain in accuracy of TRish 
with respect to TRish\_AS is within $4.6\%$. The minimum and maximum averaged final values of the 
sample size of TRish\_AS at termination were 
409 and 771, respectively.  In Figure \ref{fig4:cina}   we show the behaviour of the average training loss and testing accuracy
versus the average number of EGE for the methods equipped with their own best parameters settings.

The statistics corresponding to the best runs for both methods are reported in Table \ref{tab:cina}.
In the reported runs the sample size of TRish\_AS at termination is  approximately  4.4\% of $N$.
Concerning the step taken, consider the second row of Table 
\ref{tab:cina}; on average TRish selected the normalized step (Case 2)
$34\%$  times and the SG-type step (Case 1 and 3) $66\%$   times;  on average TRish\_AS selected 
the normalized step (Case 2)
$25\%$ times and the SG-type step (Case 3) $75\%$  times.

\begin{figure}[t!]
\begin{center}
\begin{tabular}{c}
\includegraphics[scale = 0.4]{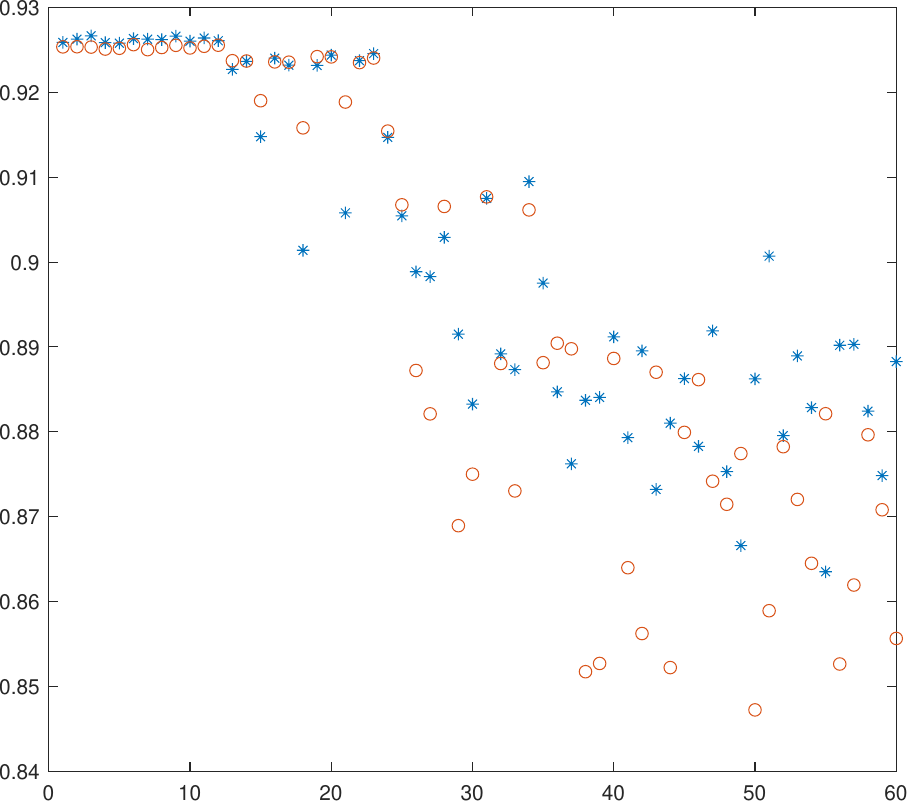}
\end{tabular}
\caption{Dataset Cina. Average testing accuracy at termination. 
TRish: symbol  ``*'', TRish\_AS: symbol ``o''. }
\label{fig1:cina}
\end{center}
\end{figure}

\begin{figure}[t!]
\begin{center}
\begin{tabular}{c}
\includegraphics[scale = 0.4]{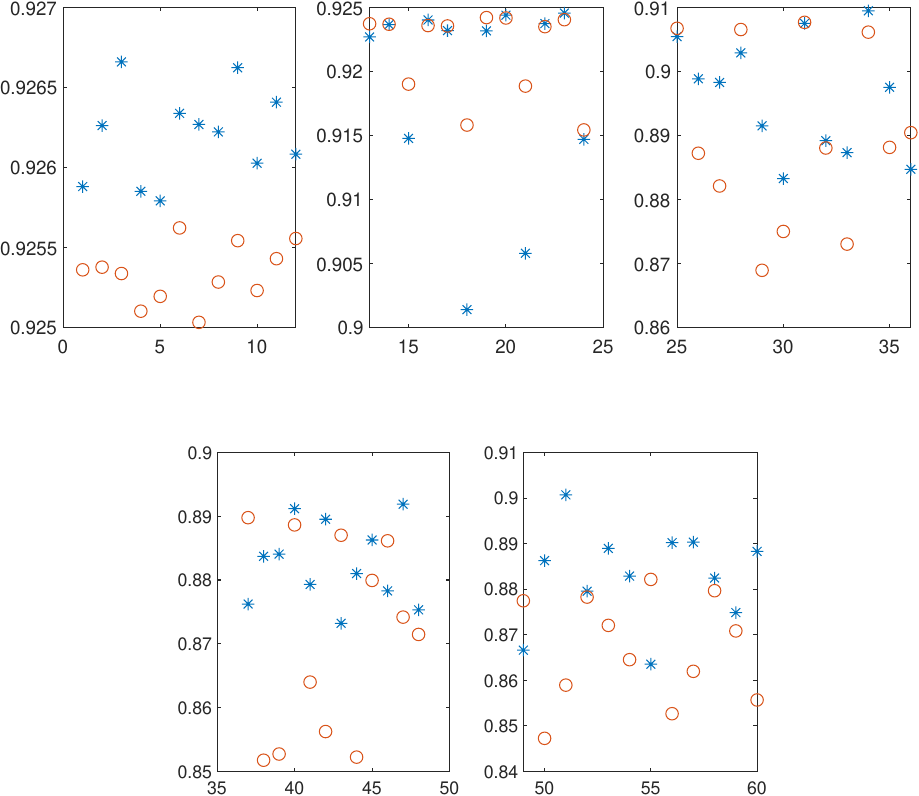}
\end{tabular}
\caption{Dataset Cina. Average testing accuracy at termination. 
Runs with fixed steplength $\alpha$ and varying $\gamma_1, \gamma_2$.
Top, from left to right:  $\alpha=10^{-1}$, $\alpha=10^{-\frac 1 2}$, $\alpha =1$. Bottom, from left to right: $\alpha=10^{\frac 1 2}$, $\alpha=10$. 
TRish: symbol  ``*'', TRish\_AS: symbol ``o''.}
\label{fig1_cinaplot4}
\end{center}
\end{figure}

\begin{figure}[t!]
\begin{center}
\begin{tabular}{c}
\includegraphics[scale = 0.4]{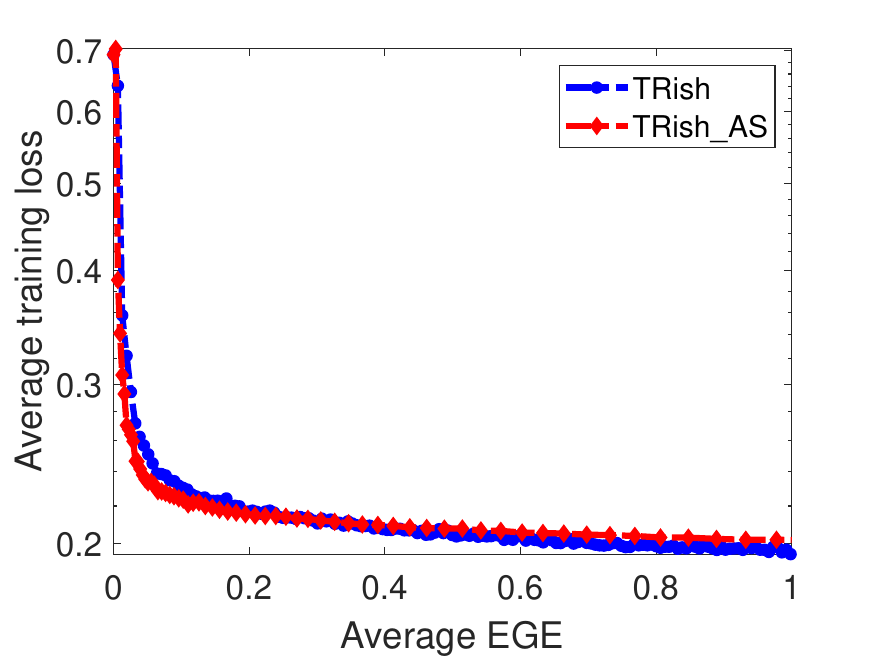}
\includegraphics[scale = 0.4]{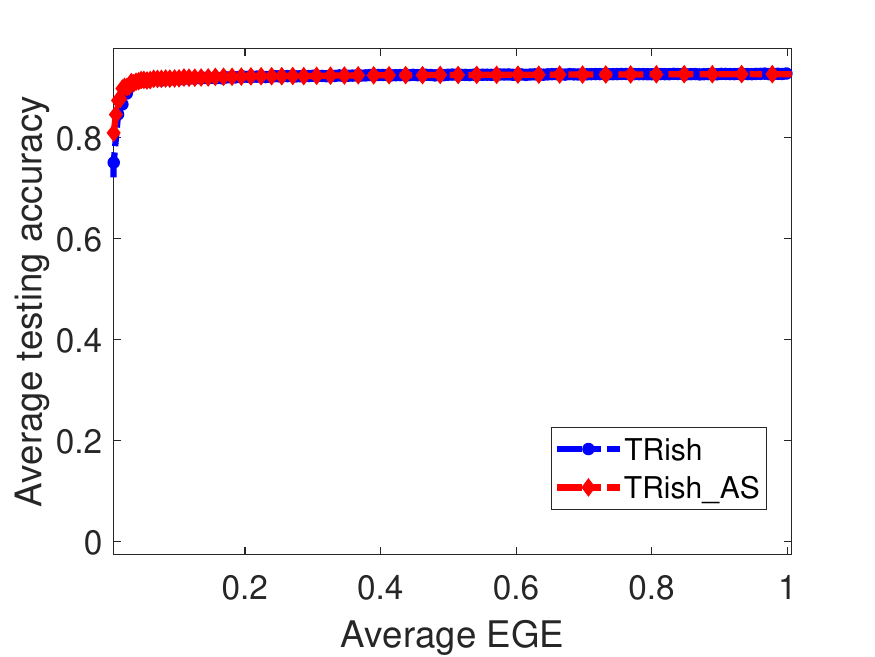}
\end{tabular}
\caption{Dataset Cina. Average training loss (left) and testing accuracy (right) versus one epoch for 
TRish  (blue) and  TRish\_AS (red).}
\label{fig4:cina}
\end{center}
\end{figure}

\begin{table}[h!]
\small
\centering
\begin{tabular}{ l|c|c|c|c}
\multicolumn{ 1}{c|}{} &\multicolumn{ 1}{c|}{$(\alpha, \gamma_1, \gamma_2)$}  & \multicolumn{ 2}{c |}{Testing accuracy} & {Average final $|S_k|$}\\ \hline
  &  & \multicolumn{ 1}{c|}{  TRish} & \multicolumn{ 1}{c|}{  TRish\_AS}  & { TRish\_AS }  \\ \hline
 Best run of TRish   & (0.1000,   18.4672,    9.2336 ) &    0.9267  &  0.9253  & 448 \\
 Best run of TRish\_AS  &   $(0.1000 ,  36.9344   , 9.2336  )$ &  0.9263  &  0.9256 & 436 
  \\
\end{tabular}  
\caption{Dataset Cina.   First row:   parameter setting corresponding to the best run for TRish in terms of testing accuracy, average testing accuracy of TRish and TRish\_AS at termination with the triplet,  average sample size of TRish\_AS at termination  with the triplet.
 Second row:   parameter setting corresponding to the best run for TRish\_AS in terms of testing accuracy, average testing accuracy of TRish and TRish\_AS at termination with the triplet,  average sample size of TRish\_AS at termination  with the triplet.}\label{tab:cina}
\end{table} 

\subsection{Neural network based classification}
In this second test, we address the problem of training a feedforward neural network for image classification. In particular, we consider the popular mnist dataset \cite{LeCun}, which consists of a training set of $N=60000$ binary images of hand-written digits from zero to nine, and a testing set of size $\bar{N} = 10000$. Our task is classifying whether an image of the testing set represents the digit ``two'' or not. We let $\{(z_i,y_i)\}_{i=1}^N$ be the training set, where $z_i\in\R^{\ell}$ is the binary image with $\ell = 784$ and $y_i\in\{0,1\}$ is the corresponding label, where $y_i=1$ detects the presence of a ``two'' and $y_i=0$ its absence. We also denote with $h(\cdot;x):\R^{\ell}\rightarrow \R$ the chosen feedforward neural network, where $x\in\R^n$ is the vector of all weights and biases. Note that the neural network has a single hidden layer equipped with $n_1=5$ nodes, whereas both the hidden and output layers are combined with the sigmoid activation function. As a result, the neural network possesses $n = (\ell+2) n_1 +1 = 3931$ parameters, which are computed by approximately solving the finite-sum minimization problem
\begin{equation}\label{eq:regression}
\color{black}\min_{x\in\R^n}F(x)\equiv -\frac 1 N\sum_{i=1}^N\left( y_i\log(h(z_i;x)) + (1-y_i)\log(1-h(z_i;x))\right).
\end{equation}
{\color{black}where the neural network is combined with the cross-entropy loss function}. We performed a similar set of experiments  as in Section \ref{sec:regression}, namely, we run both TRish and TRish\_AS for $60$ different parameter settings $(\alpha,\gamma_1,\gamma_2)$, where $\alpha\in\{0.1,10^{\frac{1}{2}},1,10^{\frac{1}{2}},10\}$, $\gamma_1\in\{\frac 4 G, \frac 8 G , \frac{16}{G}, \frac{32}{G}\}$, $\gamma_2\in\{\frac{1}{2G}, \frac{1}{G}, \frac{2}{G}\}$, being $G=0.2091$ the average norm of the stochastic gradients generated by running SG for one epoch with stepsize $\alpha=0.1$ and sample size equal to $64$. For each combination of parameters, both algorithms are run for one epoch, by using a fixed sample size $S=64$ for TRish, and an initial sample size $S_0=\min\left\{32, \left\lceil \frac{N}{100}\right\rceil\right\}$ for TRish\_AS.

As can be seen from Figure \ref{fig1:mnist} and \ref{fig1_mnistplot4}, TRish\_AS is more accurate than TRish in {\color{black} $35$ cases out of $60$; furthermore we note that the accuracy achieved by TRish\_AS is considerably higher than the accuracy obtained with TRish in several cases.  Table \ref{tab:mnist} shows the highest testing accuracy achieved when the two algorithms are equipped with their own best parameters configuration; it can be seen that TRish\_AS reached the largest testing accuracy among all the possible configurations. The minimum and maximum averaged sample sizes of TRish\_AS at termination are $77$  and
$29887$, respectively. Figure \ref{fig4:mnist} shows the decrease of the average training loss and testing accuracy versus the average number of EGE for both methods equipped with their own best parameters settings. We can observe that  the training loss in TRish\_AS shows a steady behaviour as the computational cost increases,  as opposed to the quite erratic behaviour of TRish. As for the steps taken by TRish and TRish\_AS, we report once again the statistics related to the case where both methods are equipped with the best parameters configuration for TRish\_AS (see second row of Table \ref{tab:mnist}). On average, TRish selected the normalized step (Case 2) $10\%$ times and the SG-type step (Case 1 and 3) $90\%$   times, whereas TRish\_AS selected on average
the normalized step (Case 2)
$43\%$ times and the SG-type step (Case 1 and 3) $57\%$  times.}

\begin{figure}[t!]
\begin{center}
\begin{tabular}{c}
\includegraphics[scale = 0.4]{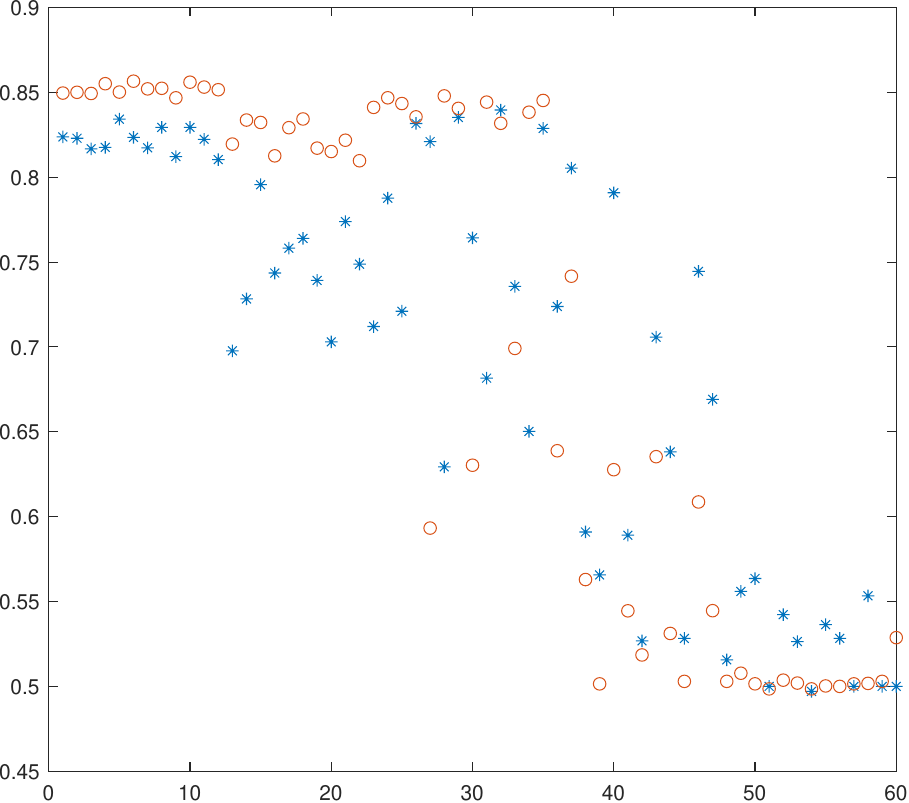}
\end{tabular}
\caption{Dataset mnist. Average testing accuracy at termination. 
TRish: symbol  ``*'', TRish\_AS: symbol ``o''. }
\label{fig1:mnist}
\end{center}
\end{figure}

\begin{figure}[t!]
\begin{center}
\begin{tabular}{c}
\includegraphics[scale = 0.4]{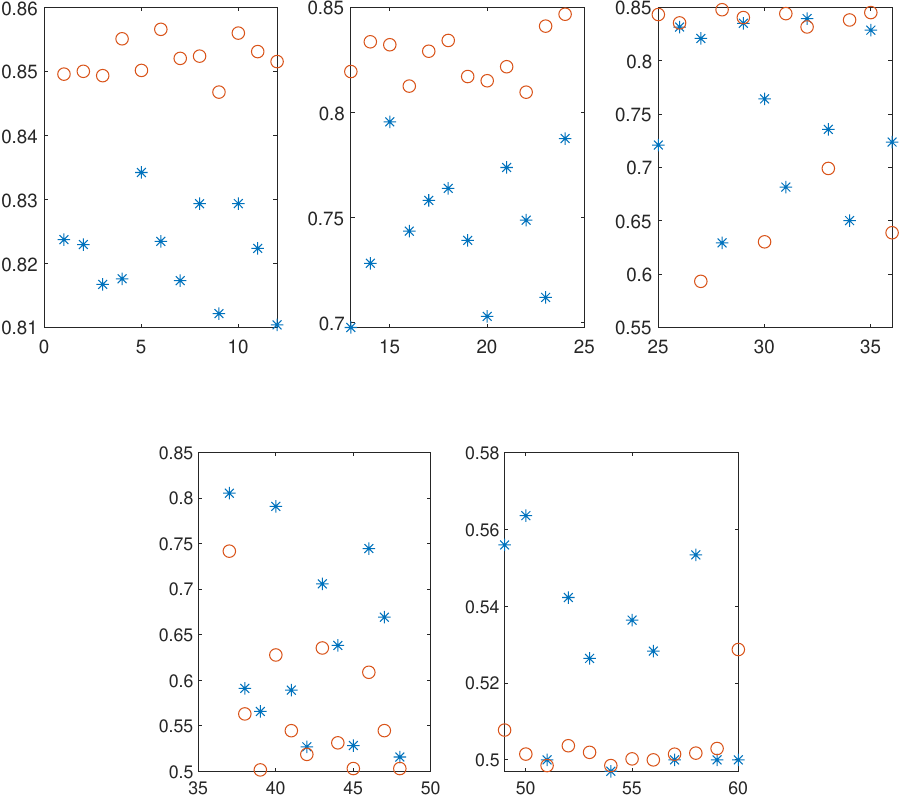}
\end{tabular}
\caption{Dataset mnist. Average testing accuracy at termination. 
Runs with fixed steplength $\alpha$ and varying $\gamma_1, \gamma_2$.
Top, from left to right:  $\alpha=10^{-1}$, $\alpha=10^{-\frac{1}{2}}$, $\alpha =1$. Bottom, from left to right: $\alpha=10^{\frac{1}{2}}$, $\alpha = 10$.
TRish: symbol  ``*'', TRish\_AS: symbol ``o''. }
\label{fig1_mnistplot4}
\end{center}
\end{figure}

\begin{figure}[t!]
\begin{center}
\begin{tabular}{c}
\includegraphics[scale = 0.4]{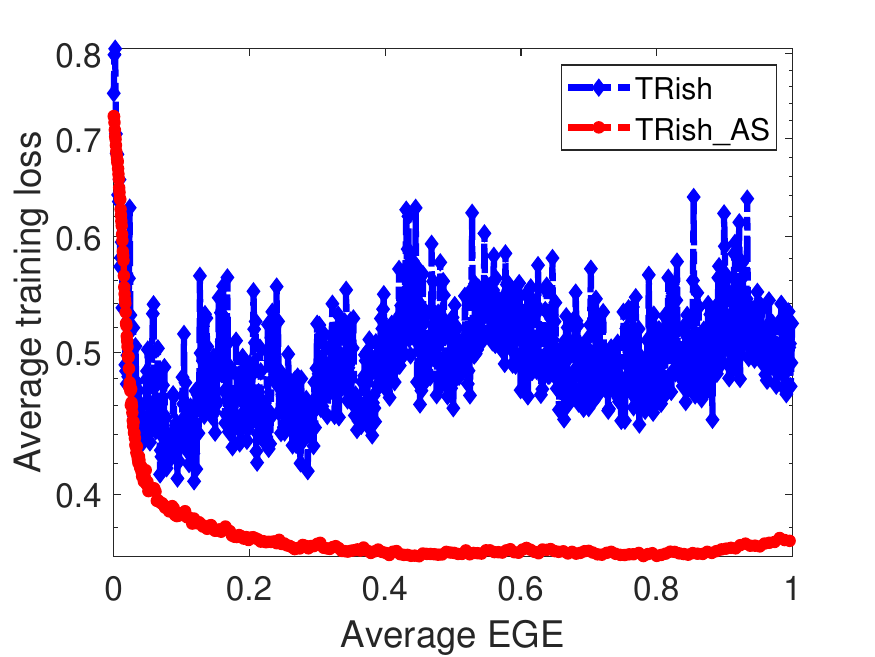}
\includegraphics[scale = 0.4]{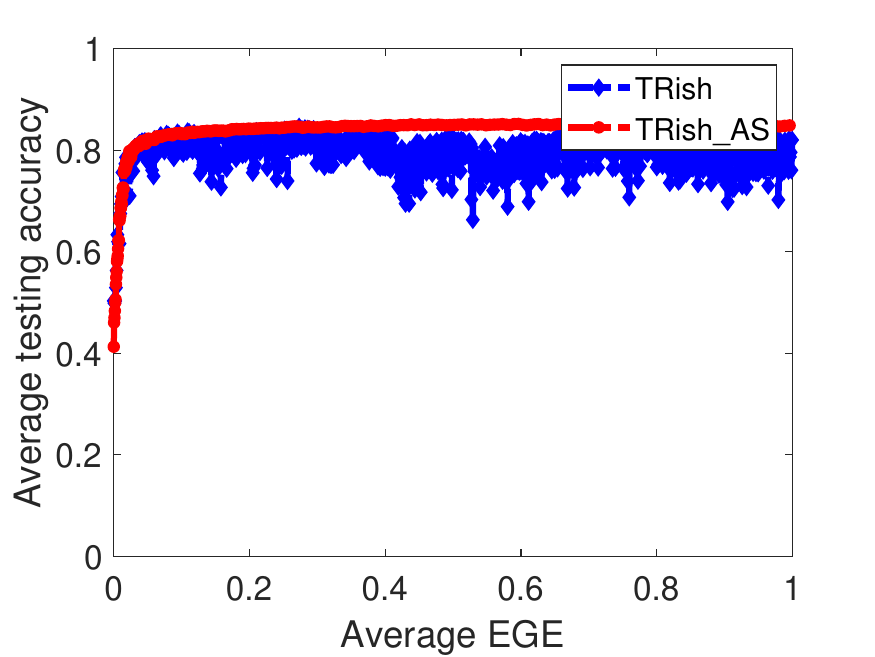}
\end{tabular}
\caption{Dataset mnist. Average training loss (left) and testing accuracy (right) versus one epoch for 
TRish (blue) and  TRish\_AS (red). }
\label{fig4:mnist}
\end{center}
\end{figure}

\begin{table}[h!]
\small
\centering
{\color{black}
\begin{tabular}{ l|c|c|c|c}
\multicolumn{ 1}{c|}{} &\multicolumn{ 1}{c|}{$(\alpha, \gamma_1, \gamma_2)$}  & \multicolumn{ 2}{c |}{Testing accuracy} & {Average final $|S_k|$}\\ \hline
  &  & \multicolumn{ 1}{c|}{  TRish} & \multicolumn{ 1}{c|}{  TRish\_AS}  & { TRish\_AS }  \\ \hline
 Best run of TRish   & $
(1.0000,    76.5184 ,  4.7824)$ &    0.8396   &    0.8318 &  1044  \\
 Best run of TRish\_AS  &   $(0.1000  ,  38.2592   ,9.5648)$ &   0.8235 &  0.8566 & 520  \\
\end{tabular} 
}
\caption{Dataset mnist.   First row:   parameter setting corresponding to the best run for TRish in terms of testing accuracy, average testing accuracy of TRish and TRish\_AS at termination with the triplet,  average sample size of TRish\_AS at termination  with the triplet.
 Second row:   parameter setting corresponding to the best run for TRish\_AS in terms of testing accuracy, average testing accuracy of TRish and TRish\_AS at termination with the triplet,  average sample size of TRish\_AS at termination  with the triplet.}  \label{tab:mnist}
\end{table}

\subsection{Neural network based regression}
We now turn our attention to an optimization problem coming from nonlinear regression.
 Let $\{(z_i,y_i)\}_{i=1}^N$ be the training set,  where $z_i\in\R^{\ell}$  represents the $i$-th feature vector, and $y_i\in\R$ is the corresponding target variable. By denoting with $h(\cdot;x):\R^{\ell}\rightarrow \R$ a prediction function, we aim at estimating the parameters vector $x\in\R^n$ by addressing an optimization problem of the form \eqref{eq:regression}.

In our tests, we use the {\sc air} dataset employed in other previous works \cite{Bellavia-et-al-22,DeVito-et-al-08}, which is available at the website \cite{UCI}. This dataset contains $9358$ instances of (hourly averaged) concentrations of polluting gases, as well as temperatures and relative/absolute air humidity levels, recorded at each hour in the period March 2004 - February 2005 from a device located in a polluted area within an Italian city. More specifically, each instance of the dataset consists of $8$ features, including the concentrations of benzene, carbon monoxide, nitrogen oxides, ozone, non-metanic hydrocarbons, nitrogen dioxide, air temperature, and relative air humidity. 

We aim at predicting the benzene concentration at a certain hour from the knowledge of the other $\ell=7$ features. After removing the instances in which the benzene concentration is missing, reducing the dataset dimension from $9357$ to $8991$, we consider $70\%$ of the dataset for training ($N=6294$), and the remaining $30\%$ for testing ($N_T=2697$). In other words, the training phase is based on the data measured in the first $9$ months, where the testing phase is performed using the data of the last  $3$ months. By letting $D=(d_{ij})\in \R^{(N+N_T)\times \ell}$ be the matrix containing all the dataset examples along its rows, and defining the values $m_j = \min\limits_{i=1,\ldots,N+N_{T}} {d_{ij}}$ and $M_j = \max\limits_{i=1,\ldots,N+N_{T}} {d_{ij}}$ for $j=1,\ldots,n$, we preprocess the data as below
$$
d_{ij} = \frac{d_{ij}-m_{j}}{M_j-m_j}, \quad i=1,\ldots, N+N_T, \ j=1,\ldots,n,
$$
so that all data values are constrained into the interval $[0,1]$.

The prediction function $h(\cdot;x)$ is chosen as a feed-forward neural network with two hidden layers of $7$ and $5$ nodes respectively, as done in \cite{DeVito-et-al-08}. Consequently, $x\in\R^{102}$ represents the vector of weights and biases of the considered neural network. Each hidden layer is equipped with the linear activation function, whereas the output layer is combined with the sigmoid activation function.

We considered the same $60$ parameter settings $(\alpha,\gamma_1,\gamma_2)$ as the ones defined in the previous sections, namely, $\alpha\in\{0.1,10^{\frac{1}{2}},1,10^{\frac{1}{2}},10\}$, $\gamma_1\in\{\frac 4 G, \frac 8 G , \frac{16}{G}, \frac{32}{G}\}$, $\gamma_2\in\{\frac{1}{2G}, \frac{1}{G}, \frac{2}{G}\}$, where $G=0.0318$ is again the average norm of stochastic gradients generated by running SG for one epoch with stepsize $\alpha=0.1$ and sample size equal to 64. For each combination of parameters, we run both TRish and TRish\_AS for one epoch, by using a fixed sample size $S=64$ for TRish, and an initial sample size $S_0=\min\left\{32, \left\lceil \frac{N}{100}\right\rceil\right\}$ for TRish\_AS.

As this is not a classification problem, we replace the testing accuracy with the evaluation of the testing loss at iteration $k$, i.e.,
\begin{equation*}
F_{S_{T}}(x_k)= \frac{1}{N_T}\sum_{i\in S_{T}} \left(y_i-h(z_i;x_k) \right)^2,
\end{equation*}   
where $N_T$ is the dimension of the testing set, and $S_T\subseteq \{1,\ldots,N\}$, $\lvert S_T \rvert = N_T$, is the set of indexes corresponding to the instances of the testing set. Note that a lower value of the testing loss corresponds to a more accurate prediction of the target variable on the testing set.  

In Figure \ref{fig1:air}, we report the average testing loss over $50$ runs at termination for each parameters setting. We can appreciate once more the good performance of TRish\_AS, which outperforms TRish  for all $60$ runs. The computational effort of TRish\_AS is overall limited, as the minimum and maximum averaged sample size at termination are $32$ and $91$, respectively. Figure \ref{fig1_airplot4} shows that TRish\_AS is completely insensitive to the choice of $\gamma_1$ and $\gamma_2$  when $\alpha=10^{-1}$, whereas TRish is always sensitive to these parameters for any choice of $\alpha$. From Table \ref{tab:air}, we see  TRish\_AS performs better than TRish even when it is equipped with the best parameters configuration for TRish.
In Figure \ref{fig4:air}   we show the behaviour of the average training and testing loss versus the average number
of effective gradient evaluations (EGE) for the methods equipped with their own best parameters settings.  We can appreciate an accelerated behaviour of TRish\_AS with respect to TRish in both training and testing loss. The statistics corresponding to the best runs for both methods are reported in Table \ref{tab:air}.
Concerning the step taken, consider the second row of Table 
\ref{tab:cina}; on average TRish selected the normalized step (Case 2)
$81\%$  times and the SG-type step (Case 3) $19\%$   times;  on average TRish\_AS selected 
the normalized step (Case 2)
$77\%$ times and the SG-type step (Case 1 and 3) $23\%$  times.
  
\begin{figure}[t!]
\begin{center}
\begin{tabular}{c}
\includegraphics[scale = 0.4]{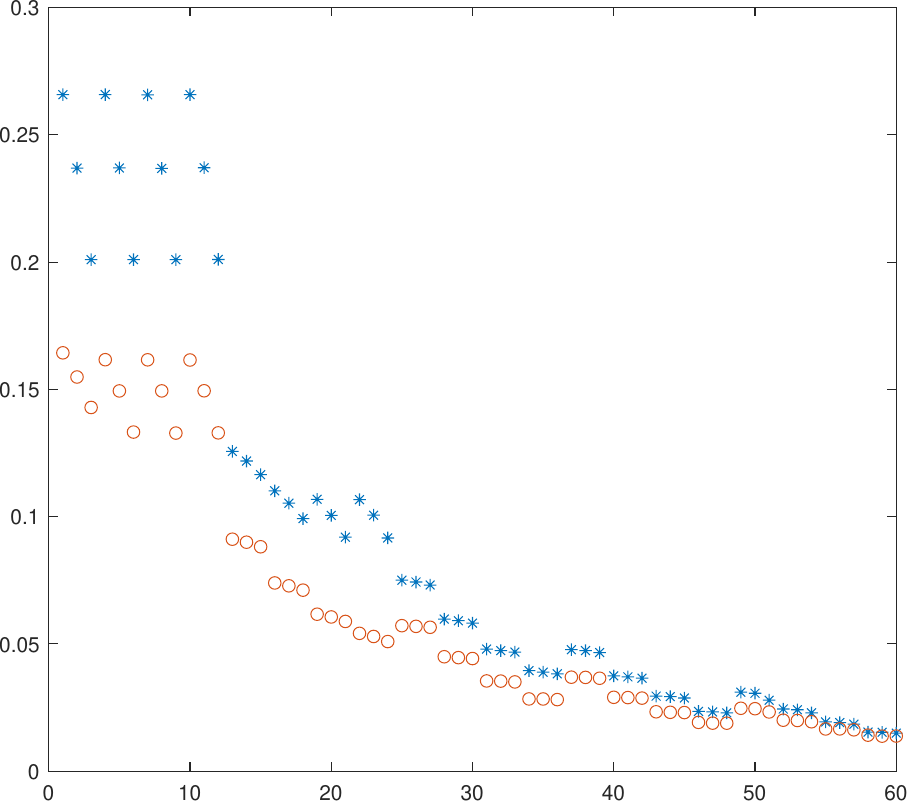}
\end{tabular}
\caption{Dataset {\sc air}. Average testing loss at termination. 
TRish: symbol  ``*'', TRish\_AS: symbol ``o''. }
\label{fig1:air}
\end{center}
\end{figure}

\begin{figure}[t!]
\begin{center}
\begin{tabular}{c}
\includegraphics[scale = 0.4]{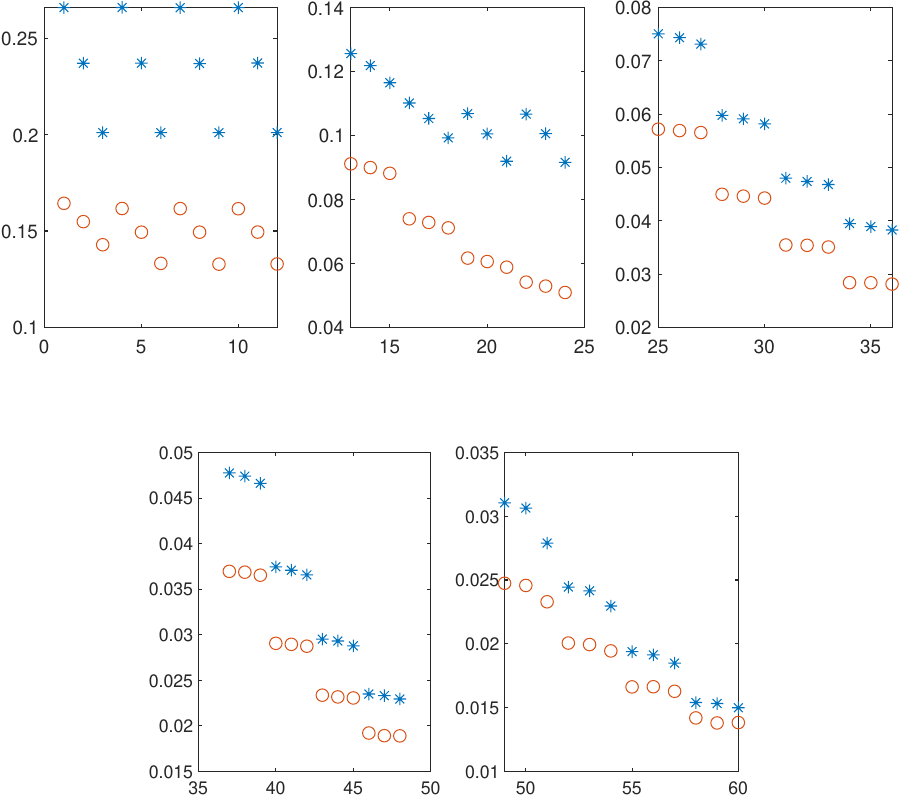}
\end{tabular}
\caption{Dataset {\sc air}. Average testing loss at termination. 
Runs with fixed steplength $\alpha$ and varying $\gamma_1, \gamma_2$.
Top, from left to right:  $\alpha=10^{-1}$, $\alpha=10^{-\frac{1}{2}}$, $\alpha=1$. Bottom, from left to right: $\alpha=10^{\frac{1}{2}}$, $\alpha = 10$. TRish: symbol  ``*'', TRish\_AS: symbol ``o''. }
\label{fig1_airplot4}
\end{center}
\end{figure}

\begin{figure}[t!]
\begin{center}
\begin{tabular}{c}
\includegraphics[scale = 0.4]{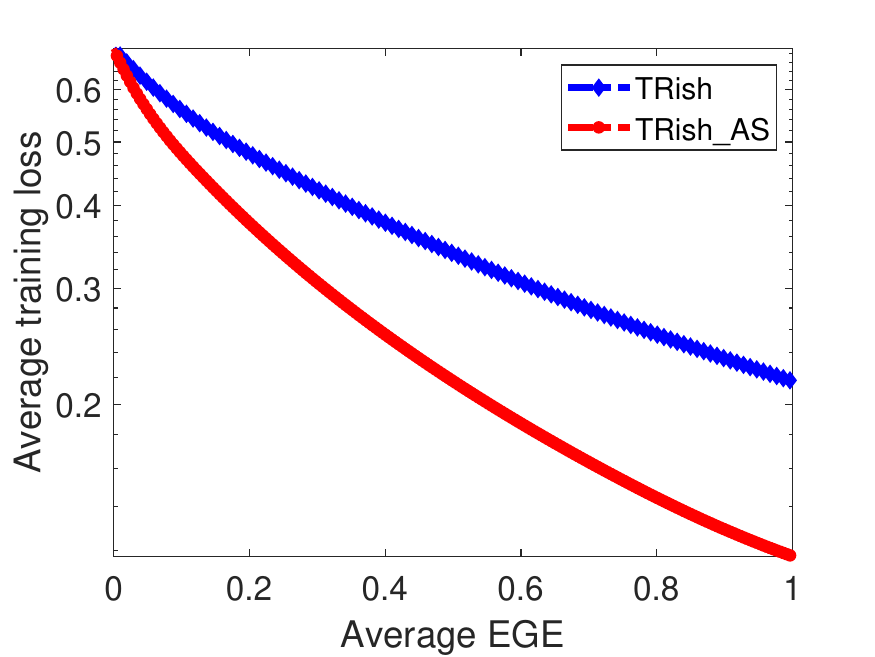}
\includegraphics[scale = 0.4]{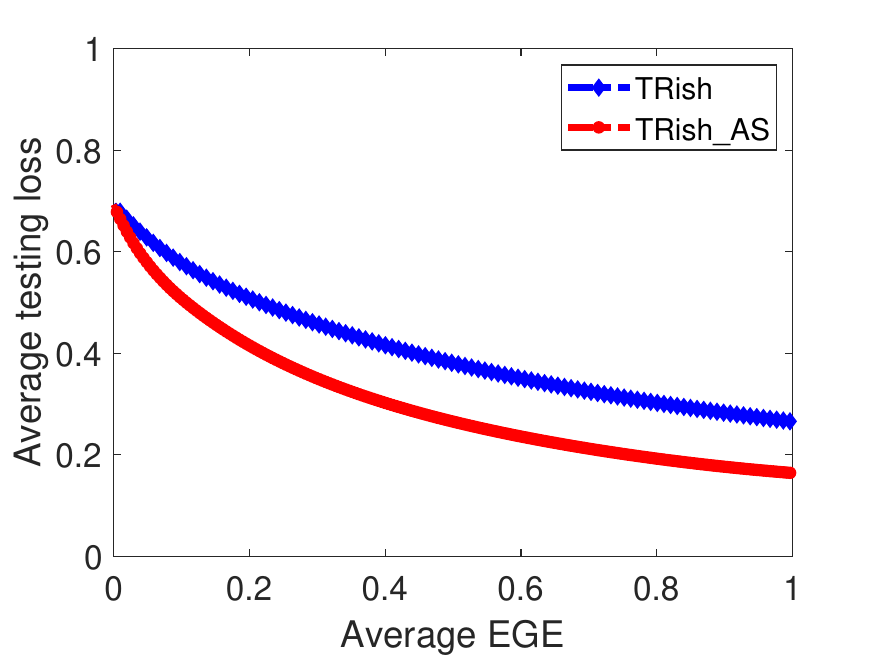}
\end{tabular}
\caption{ Dataset air. Average training loss (left) and testing loss (right) versus one epoch for 
TRish (blue) and  TRish\_AS (red). }
\label{fig4:air}
\end{center}
\end{figure}

\begin{table}[h!]
\small
\centering
\begin{tabular}{ l|c|c|c|c}
\multicolumn{ 1}{c|}{} &\multicolumn{ 1}{c|}{$(\alpha, \gamma_1, \gamma_2)$}  & \multicolumn{ 2}{c |}{Testing loss} & {Average final $|S_k|$}\\ \hline
  &  & \multicolumn{ 1}{c|}{  TRish} & \multicolumn{ 1}{c|}{  TRish\_AS}  & { TRish\_AS }  \\ \hline
 Best run of TRish   & $(10,  1005.4492, 62.8406)$ &    $0.01498$  &   $0.01382$ & $91$ \\
 Best run of TRish\_AS  &   $(10,  1005.4492, 31.4203)$ &   $0.01531$ &   $0.01379$   & $83$  \\
\end{tabular}  
\caption{Dataset {\sc air}.   First row:   parameter setting corresponding to the best run for TRish in terms of testing loss, average testing loss of TRish and TRish\_AS at termination with the triplet,  average sample size of TRish\_AS at termination  with the triplet.
 Second row:   parameter setting corresponding to the best run for TRish\_AS in terms of testing loss, average testing loss of TRish and TRish\_AS at termination with the triplet,  average sample size of TRish\_AS at termination  with the triplet.}\label{tab:air}
\end{table}

\section{Conclusions}\label{sec:conclusions}
In this work a theoretical and experimental analysis of the stochastic trust-region based TRish algorithm was carried out.
Our analysis complements and partly improves upon the existing results.   Our implementation is based on a dynamic choice
of the sample size, and provides good results with respect to using a prefixed and constant sample size. In particular, the performance of the algorithm is less sensitive to the choice of the parameters governing the form of the step, once the stochastic gradient is formed.

\end{document}